\documentclass[a4paper,times,3p]{elsarticle}
\usepackage[utf8]{inputenc}

\usepackage[T1]{fontenc}

\usepackage{microtype}

\usepackage{booktabs}

\usepackage{amsmath}
\usepackage{amssymb}
\usepackage{amsthm}
\usepackage{bm}
\usepackage{graphicx}
\usepackage{graphbox}
\usepackage{subfig}
\usepackage{epstopdf}
\usepackage{psfrag}
\usepackage{color}
\usepackage{hyperref}
\usepackage{diagbox}
\usepackage{enumitem}
\usepackage{dirtytalk}
\usepackage{tikz}
\usetikzlibrary{shapes, arrows}

\usepackage{mathtools}

\definecolor{cmykcyan}{cmyk}{1,0,0,0}
\definecolor{cmykred}{cmyk}{0,1,1,0}
\definecolor{cmykblack}{cmyk}{0,0,0,1}

\newtheorem*{remark}{Remark}

\usepackage[font=footnotesize,labelfont=bf]{caption}
\usepackage{multirow}

\usepackage{pxfonts}

\usepackage{import}

\graphicspath{./graphics/}

\usepackage[textsize=scriptsize]{todonotes}
\usepackage{setspace}

\makeatletter
\def\ps@pprintTitle{%
 \let\@oddhead\@empty
 \let\@evenhead\@empty
 \def\@oddfoot{}%
 \let\@evenfoot\@oddfoot}
\makeatother

\newcommand{\balpha}{\ensuremath{\boldsymbol{\alpha}}}
\newcommand{\dA}{\ensuremath{\mathbf{d}_{\mathcal{A}}}}
\newcommand{\cA}{\ensuremath{\mathbf{c}_{\mathcal{A}}}}
\newcommand{\cB}{\ensuremath{\mathbf{c}_{\mathcal{B}}}}
\newcommand{\cI}{\ensuremath{\mathbf{c}_{\mathcal{I}}}}
\DeclareMathAlphabet{\mymathbb}{U}{BOONDOX-ds}{m}{n}

\begin{document}

\tikzstyle{block} = [draw, fill=blue!20, rectangle, 
    minimum height=3em, minimum width=6em]
\tikzstyle{sum} = [draw, fill=blue!20, circle, node distance=1cm]
\tikzstyle{input} = [coordinate]
\tikzstyle{output} = [coordinate]
\tikzstyle{pinstyle} = [pin edge={to-,thin,black}]

\begin{frontmatter}

\title{The Role of PDE-Based Parameterization Techniques in Gradient-Based IGA Shape Optimization Applications}

\author[add1]{Jochen Hinz\corref{cor1}}
\ead{j.p.hinz@tudelft.nl}
\author[add2]{Andrzej Jaeschke\corref{}}
\ead{andrzej.jaeschke@p.lodz.pl}
\author[add1]{Matthias M\"oller\corref{}}
\ead{m.moller@tudelft.nl}
\author[add1]{Cornelis Vuik \corref{}}
\ead{C.Vuik@tudelft.nl}

\cortext[cor1]{Corresponding author}

\address[add1]{Department of Applied Mathematics\\
	Delft University of Technology, 2628 XE Delft, Netherlands}
\address[add2]{Institute of Turbomachinery \\ Lodz University of Technology, Wolczanska 219/223, 90-924 Lodz, Poland}

\begin{abstract}
This paper proposes a shape optimization algorithm based on the principles of \textit{Isogeometric Analysis} (IGA) in which the parameterization of the geometry enters the problem formulation as an additional PDE-constraint. Inspired by the \textit{isoparametric principle} of IGA, the parameterization and the governing state equation are treated using the same numerical technique. This leads to a scheme that is comparatively easy to differentiate, allowing for a fully symbolic derivation of the gradient and subsequent gradient-based optimization. To improve the efficiency and robustness of the scheme, the basis is re-selected during each optimization iteration and adjusted to the current needs. The scheme is validated in two test cases.
\end{abstract}

\end{frontmatter}

\section{Introduction}
\label{sect:Introduction}

Isogeometric analysis (IGA) was introduced by Hughes et al. in \cite{hughes2005isogeometric} as a numerical technique that bridges the gap between computer aided design (CAD) and the numerical analysis of Partial Differential Equations (PDEs).
This is accomplished by using the same function space to represent the geometry $\Omega$ and to discretize the PDE-problem posed over $\Omega$. \\ 
Most of the available CAD software generates no more than a spline-based description of the boundary contours $\partial \Omega$ of $\Omega$. Therefore, suitable parametrization algorithms are indispensable for generating bijective (folding-free), analysis-suitable geometry parameterizations from the boundary CAD data. \\
The parameteric quality of the mapping has a profound impact on the numerical accuracy of the isogeometric analysis \cite{xu2010optimal}. Therefore, besides bijectivity, proficient parameterization algorithms aim at generating parameterizations of high numerical quality. \\
A variety of parametrization techniques have been proposed in the literature such as Coon's Patch \cite{farin1999discrete}, Linear Spring \cite{gravesen2012planar} and approaches based on (constrained and unconstrained) quality cost function optimization \cite{gravesen2012planar, xu2011parameterization, falini2015planar}. While mappings based on Coon's Patch and Linear Spring follow from a closed-form expression and are hence cheap to compute and straightforwardly differentiable, they often lead to folded (non-bijective) mappings. The same is true for unconstrained optimization. Constrained optimization approaches on the other hand typically have a higher success rate. However, this comes at the expense of a large number of (constrained) iterations (typically about $\sim 30$) and the notorieties associated with nonconvex optimization, such as the danger of getting stuck in local minima. A third class of approaches attempts to generate a mapping whose inverse is composed of harmonic functions on $\Omega$. This approach is based on the observation that harmonic functions exhibit a large degree of smoothness, which benefits the numerical quality of the resulting mapping. Furthermore, it can be shown that inversely harmonic mappings (IHMs) are bijective, thanks to the maximum principle \cite{rado1926aufgabe, kneser1926losung}. Many approaches for approximating IHMs have been proposed in the literature \cite{nguyen2010parameterization, falini2015planar}, notably the PDE-based approach called Elliptic Grid Generation (EGG) \cite{azarenok2009generation, hinz2018elliptic}. EGG is of particular interest in shape optimization problems thanks to the parametric smoothness and bijectivity of IHMs as well as differentiability, made possible by the PDE-based problem formulation. \\
Traditionally, IGA parametrizations are built from tensor-product spline spaces. Unfortunately, structured spline technologies do not allow for local refinement. This may result in infeasibly-large function spaces. Therefore, unstructured spline technologies such as THB-splines \cite{giannelli2012thb} are gaining an increased amount interest in the IGA community, thanks to local refinement. An EGG-based planar parametrization framework that supports THB-splines has been proposed in \cite{hinz2020goaloriented}. \\
Since its birth in $2005$, IGA was successfully applied to wide variety of problems including: thermal analysis \cite{kacprzyk2014isogeometric}, linear elasticity problems \cite{hughes2005isogeometric}, structural vibrations \cite{cottrell2006isogeometric}, incompressible flows \cite{bazilevs2008nurbs} and inviscid compressible flows \cite{Jaeschke2015}. As a mature numerical method, it is ready to be used in more complex industrial processes. As a result, several publications that apply IGA to shape optimization problems have appeared in the literature \cite{yoon2013isogeometric, nortoft2013isogeometric, wall2008isogeometric, manh2011isogeometric}. Combining IGA and shape optimization is very appealing as the spline-based description of $\partial \Omega$ can be used directly to compute a mapping for $\Omega$, completely bypassing the need to first convert $\partial \Omega$ into a piecewise-linear curve that acts as input for classical mesh generators. \\
There are two main groups of shape optimization algorithms: gradient-free (like for example genetic algorithms \cite{holland1984genetic}) and gradient-based methods (for example interior point methods \cite{wachter2006implementation, biegler2009large}). The latter group generally requires fewer underlying PDE evaluations at the expense of having to compute the gradient of the objective function during each iteration. Therefore, differentiability of the IGA parametrization algorithm constitutes a significant advantage. An additional feature of differentiability is efficiency: as the inner control points are a smooth function of the boundary control points, there is no need for full remeshing after each iteration since cheaper mesh update strategies can be employed. This is also true for settings in which the boundary contours change as a smooth function of time. \\
In order to combine the appealing features of EGG and THB-enabled local refinement, this paper adopts the parameterization framework proposed in \cite{hinz2020goaloriented} and presents an IGA-based shape optimization algorithm in which the parameterization is added to the optimization problem formulation in the form of an additional PDE-constraint. In line with the \textit{isoparametric principle} of IGA, we numerically treat this additional constraint in the same way as the governing quantity (temperature, pressure, etc) of the underlying optimization problem. Including the mapping explicitly as a PDE-constraint facilitates differentiation, allowing for gradient-based optimization, while also guaranteeing analysis-suitability, thanks to the bijectivity of IHMs. To improve the efficiency, the proposed algorithm employs THB-enabled adaptive local refinement strategies during every optimization iteration, resulting in a variable discretization basis. We validate the proposed methodology by presenting two test cases.
\section{Notation}
In this work, we denote vectors in boldface while matrices receive a capital letter and may furthermore be enclosed in square brackets for better readability. The $i$-th entry of vector $\mathbf{x}$ is denoted by $\mathbf{x}_i$ or simply $x_i$ and similarly for the $ij$-th entry of matrices. We make extensive use of vector derivatives. Here, we interchangeably use the denotation
\begin{align}
    \left[ \partial_{\mathbf{t}} \mathbf{x} \right] \equiv \left[ \frac{\partial \mathbf{x}}{\partial \mathbf{t}} \right], \quad \text{with} \quad \left[ \frac{\partial \mathbf{x}}{\partial \mathbf{t}} \right]_{ij} = \frac{\partial x_i}{\partial t_j}
\end{align}
for the partial derivative and similarly for the total derivative. In the case of taking the derivative of a scalar, brackets are avoided. However, the argument is treated as a $1 \times 1$ matrix and hence the derivative has dimension $(1, m)$, where $m$ is the dimension of $\mathbf{t}$.
\section{Problem Formulation}
\label{sect:Problem_Formulation}
We are considering the shape optimization problem of a planar domain $\Omega(\balpha)$ whose contours $\partial \Omega(\balpha)$ are parameterized by the $n$-tuple of design variables $\balpha = (\alpha_1, \ldots, \alpha_n)$. If the design variables are taken from the design space $\boldsymbol{\mathcal{\lambda}}$, the optimization problem reads:
\begin{align}
\label{eq:Shape_Optimization}
\begin{split}
    J\left(u^{\balpha}, \Omega^{\balpha}, \balpha \right) & \rightarrow \min_{\balpha} \\
    \text{s.t.} \quad g_i(u^{\balpha}, \Omega^{\balpha}, \balpha) & \geq 0, \quad \forall i \in \{1, \ldots, N_{\neq} \} \\
                \quad h_j(u^{\balpha}, \Omega^{\balpha}, \balpha) & = 0, \quad \forall j \in \{1, \ldots, N_{=} \} \\
    \balpha & \in \boldsymbol{\mathcal{\lambda}},
\end{split}
\end{align}
where the $g_i$ and $h_j$ are problem-specific constraints. Here, $J(\cdot, \cdot, \cdot)$ denotes the objective function and $u^{\balpha}: \Omega \rightarrow \mathbb{R}$ some state variable whose physical meaning depends on the application (temperature, pressure, etc). We regard $u^{\balpha}$ as a scalar quantity for convenience. However, generalizations to vectorial quantities are straightforward. Note that the dependencies of the variables contained in $J(\cdot, \cdot, \cdot)$ are concatenated in descending order, i.e., in general $u^{\balpha} = u^{\balpha}(\Omega^{\balpha}(\balpha), \balpha)$ and $\Omega^{\balpha} = \Omega^{\balpha}(\balpha)$. The state variable $u^{\balpha}$ follows from a PDE-problem posed over $\Omega^{\balpha}$ and may contain additional dependencies on $\balpha$ (such as source terms), hence the dependency on the tuple $(\Omega^{\balpha}, \balpha)$. Tackling (\ref{eq:Shape_Optimization}) computationally requires introducing a bijective geometry parameterization $\mathbf{x}^{\balpha}: \hat{\Omega} \rightarrow \Omega^{\balpha}$, where $\hat{\Omega}$ denotes the computational domain which is assumed to be static. Here, we restrict ourselves to geometries that are topologically equivalent to $\hat{\Omega} = [0, 1]^2$ for convenience. However the generalization to multipatch settings is straightforward. Let
\begin{align}
\label{eq:Dirichlet_basis_u}
    \mathcal{U}^{f} = \{v \in \mathcal{U} \enskip \vert \enskip v = f \text{ on } \partial \Omega^{\boldsymbol{\alpha}}_D \}
\end{align}
for some suitably-chosen test space $\mathcal{U}$ and some $\partial \Omega_D^{\balpha} \subseteq \partial \Omega^{\balpha}$ on which Dirichlet data is prescribed. Deriving the weak form of the PDE-problem governing $u^{\balpha}$ leads to
\begin{align}
\label{eq:u_weak_not_discretized}
    \text{find } u^{\balpha} \in \mathcal{U}^{u_D^{\balpha}} \quad \text{s.t.} \quad B \left(u^{\balpha},\mathbf{x}^{\balpha}, \balpha, \phi \right) = 0, \quad \forall \phi \in \mathcal{U}^{0},
\end{align}
for some differential form $B(\cdot, \cdot, \cdot, \cdot)$. Here, $u_D^{\balpha}$ denotes the Dirichlet data as a function of the design variables. By introducing the mapping $\mathbf{x}^{\balpha}$, the objective function takes the form
\begin{align}
\label{eq:J_of_x}
    J(u^{\balpha}, \Omega^{\balpha}, \balpha) \rightarrow J\left(u^{\balpha}, \mathbf{x}^{\balpha}, \balpha \right) \equiv J^{\balpha},
\end{align}
where $u^{\balpha}$ satisfies (\ref{eq:u_weak_not_discretized}). With the dependencies of $u^{\balpha}$ and $\mathbf{x}^{\balpha}$ in mind, the gradient of (\ref{eq:J_of_x}) reads:
\begin{align}
\label{eq:grad_J_x}
    \frac{\mathrm{d} J^{\balpha}}{\mathrm{d} \balpha} & = \frac{\partial J^{\balpha}}{\partial u^{\balpha}} \left( \frac{\partial u^{\balpha}}{\partial \mathbf{x}^{\balpha}} \frac{\mathrm{d} \mathbf{x}^{\balpha}}{\mathrm{d} \balpha} + \frac{\partial u^{\balpha}}{\partial \balpha} \right) + \frac{\partial J^{\balpha}}{\partial \mathbf{x}^{\balpha}} \frac{\mathrm{d} \mathbf{x}^{\balpha}}{\mathrm{d} \balpha} + \frac{\partial J^{\balpha}}{\partial \balpha}.
\end{align}
We see that (\ref{eq:grad_J_x}) requires taking the derivative of $\mathbf{x}^{\balpha}$ with respect to $\balpha$, while the state variable $u^{\balpha}$ needs to be differentiable with respect to $\mathbf{x}^{\balpha}$. These two derivatives often constitute the most challenging step in computing the gradient because differentiating $\mathbf{x}^{\balpha}$ or with respect to $\mathbf{x}^{\balpha}$ can be nontrivial, depending on the parameterization technique used. On the other hand, differentiation with respect to $u^{\balpha}$ is relatively straightforward because the implicit function theorem can be used on (\ref{eq:u_weak_not_discretized}). Hence, if we take $\mathbf{x}^{\balpha}$ as the solution of a PDE problem, differentiation is simplified, allowing for a symbolic derivation of all terms involved in (\ref{eq:grad_J_x}). To this end, we adopt the principles of \textit{Elliptic Grid Generation}, which will be the topic of the next section.
\section{Elliptic Grid Generation}
\label{sect:EGG}
Elliptic grid generation (EGG) is a PDE-based technique aimed at generating analysis-suitable geometry parameterizations $\mathbf{x}^{\balpha}: \hat{\Omega} \rightarrow \Omega^{\balpha}$ given only a parametric description of the boundary contours $\partial \Omega^{\balpha}$ as a function of the state vector $\balpha$. Let the free topological variables in $\hat{\Omega}$ be given by the tuple $\boldsymbol{\xi} = (\xi_1, \xi_2) = (\xi, \eta)$. Then, the equations of EGG read \cite{hinz2020goaloriented}:
\begin{align}
\label{eq:inverse_mapping_approach_PDE_scaled}
    A(\mathbf{x}^{\balpha}) \colon H(\mathbf{x}_i^{\balpha}) & = 0 \quad \text{in } \hat{\Omega}, \quad \text{for } i \in \{1, 2\} \quad \text{s.t.} \quad \mathbf{x}^{\balpha} \vert_{\partial \hat{\Omega}} = \partial \Omega^{\balpha},
\end{align}
where
\begin{align}
\label{eq:A}
H(u)_{ij} \equiv \frac{\partial^2 u}{\partial \xi_i \partial \xi_j} \quad \text{and} \quad A(\mathbf{x}^{\balpha}) = \frac{1}{g_{11} + g_{22} + \epsilon} \begin{pmatrix} g_{22} & - g_{12} \\
                                                                         - g_{12} & g_{11} \end{pmatrix},
\end{align}
with $g_{ij} = \mathbf{x}_{\xi_i}^{\balpha} \cdot \mathbf{x}_{\xi_j}^{\balpha}$ the entries of the metric tensor and $\epsilon$ a small positive constant (typically, we take $\epsilon = 10^{-4}$). Here, $A \colon B$ denotes the Frobenius inner product between matrices $A$ and $B$. The solution of (\ref{eq:inverse_mapping_approach_PDE_scaled}) is a mapping $\mathbf{x}^{\balpha}$ whose inverse $(\mathbf{x}^{\balpha})^{-1}$ constitutes a pair of harmonic functions on $\Omega^{\balpha}$. As $(\mathbf{x}^{\balpha})^{-1}$ maps into a convex computational domain $\hat{\Omega}$, it follows from the maximum principle that $\mathbf{x}^{\balpha}$ is a bijection between $\hat{\Omega}$ and $\Omega^{\balpha}$, where $\mathbf{x}^{\balpha} \vert_{\partial \hat{\Omega}}$ parameterizes $\partial \Omega^{\balpha}$ \cite{rado1926aufgabe, kneser1926losung}. This property justifies limiting the choice of $\mathbf{x}^{\balpha}$ from the set of bijective parameterizations to the subset of mappings that satisfy (\ref{eq:inverse_mapping_approach_PDE_scaled}). \\
While many alternative approaches that do not require solving a PDE problem have been proposed in the literature \cite{nguyen2010parameterization, falini2015planar}, we base a computational approach on (\ref{eq:inverse_mapping_approach_PDE_scaled}) since it facilitates differentiating $\mathbf{x}^{\balpha}$ with respect to the design variables (see (\ref{eq:grad_J_x})). \\
For a viable computational approach, we derive the weak counterpart of (\ref{eq:inverse_mapping_approach_PDE_scaled}). Here, we adopt the approach from \cite{hinz2018elliptic}. Given a differential function $\mathbf{d}^{\balpha}: \partial \hat{\Omega} \rightarrow \mathbb{R}^2$ that parameterizes $\partial \Omega^{\balpha}$, the mapping $\mathbf{x}^{\balpha}$ is the solution of:
\begin{align}
\label{eq:x_weak_not_discretized}
    \text{Find } \mathbf{x}^{\balpha} \in \boldsymbol{\mathcal{V}}^{\mathbf{d}^{\balpha}} \quad \text{s.t.} \quad F(\mathbf{x}^{\balpha}, \boldsymbol{\sigma}) = 0, \quad \forall \boldsymbol{\sigma} \in \boldsymbol{\mathcal{V}}^{\mathbf{0}},
\end{align}
with
\begin{align}
\label{eq:F_operator}
    F(\mathbf{x}^{\balpha}, \boldsymbol{\sigma}) = \int_{\hat{\Omega}} \boldsymbol{\sigma}_i A(\mathbf{x}^{\balpha}) \colon H(\mathbf{x}_i^{\balpha}) \mathrm{d} S,
\end{align}
where we have made use of the Einstein summation convention. Furthermore, in (\ref{eq:x_weak_not_discretized}) we used
\begin{align}
\label{eq:test_space_x_f}
    \boldsymbol{\mathcal{V}}^{\mathbf{f}} \equiv \{ \mathbf{v} \in \mathcal{V}^2 \enskip \vert \enskip \mathbf{v} = \mathbf{f} \text{ on } \partial \hat{\Omega} \}, \quad \text{with} \quad \mathcal{V} = H^2(\hat{\Omega}) \quad \text{and} \quad \mathbf{f} \in \mathcal{V}^2.
\end{align}
The discretization of (\ref{eq:x_weak_not_discretized}) follows straightforwardly from replacing $\mathcal{V}$ by the finite-dimensional $\mathcal{V}_h^{\balpha} \subset \mathcal{V}$ in (\ref{eq:test_space_x_f}). We denote the resulting space by $\boldsymbol{\mathcal{V}}^{\balpha, \mathbf{f}}_h$. As $\mathbf{f} \in \mathcal{V}_h^{\balpha} \times \mathcal{V}_h^{\balpha} \equiv \boldsymbol{\mathcal{V}}_h^{\balpha}$ by assumption, the discretization additionally requires replacing the Dirichlet data $\mathbf{d}^{\balpha}$ by a proper collocation $\mathbf{d}^{\balpha}_h \in \boldsymbol{\mathcal{V}}_h^{\balpha}$. As such, the fully discretized problem reads:
\begin{align}
\label{eq:x_weak_discretized}
    \text{Find } \mathbf{x}^{\balpha}_h \in \boldsymbol{\mathcal{V}}^{\balpha, \mathbf{d}^{\balpha}_h}_h \quad \text{s.t.} \quad F(\mathbf{x}^{\balpha}_h, \boldsymbol{\sigma}_h) = 0, \quad \forall \boldsymbol{\sigma}_h \in \boldsymbol{\mathcal{V}}^{\balpha, \mathbf{0}}_h.
\end{align}
\begin{remark}
Due to the appearance of second order derivatives in (\ref{eq:x_weak_discretized}), we have to assume that $\mathbf{x}^{\balpha}_h$ is built from a space with global $C^1$-continuity. For an approach that allows for lower regularity (and is hence compatible with multipatch parameterizations), we refer to \cite{hinz2019iga}.
\end{remark}
Since (\ref{eq:x_weak_discretized}) is a nonlinear root-finding problem, we tackle it with a Newton-based iterative approach. Unlike $\mathbf{x}^{\balpha}$, its discretized counterpart $\mathbf{x}^{\balpha}_h$ may fold due to the truncation error introduced by the numerical scheme. Grid folding can be repaired by refining $\mathcal{V}_h^{\balpha}$ in the affected regions and recomputing $\mathbf{x}_h^{\balpha}$ from the enriched space. This makes using an unstructured spline technology like THB-splines particularly appealing, thanks to local refinement. For more details on the choice of $\mathcal{V}_h^{\balpha}$ and the computational approach, we refer to \cite{hinz2020goaloriented}.
\section{Computational Approach}
\label{sect:Computational_Approach}
In this section we propose a computational approach for numerically treating the optimization problem (\ref{eq:Shape_Optimization}).

\subsection{Discretization}
\label{subsect:Discretization}
We discretize the optimization problem (\ref{eq:Shape_Optimization}) by approximating
\begin{align}
\label{J_discrete}
    J(u^{\balpha}, \mathbf{x}^{\balpha}, \balpha) \simeq J(u^{\balpha}_h, \mathbf{x}^{\balpha}_h, \balpha) \equiv J_h^{\balpha},
\end{align}
where $u^{\balpha}_h \in \mathcal{U}_h^{\balpha}$ is the solution of the discretized weak state equation (\ref{eq:u_weak_not_discretized}) while $\mathbf{x}^{\balpha}_h \in \boldsymbol{\mathcal{V}}^{\balpha}_h$ is the solution of (\ref{eq:x_weak_discretized}) for given $\balpha$. Here $\Omega^{\balpha}_h$ is parameterized by $\mathbf{x}^{\balpha}_h$ and approximates the domain $\Omega^{\balpha}$ whose contours are parameterized by the $\balpha$-differentiable $\mathbf{d}^{\balpha}: \partial \hat{\Omega} \rightarrow \mathbb{R}^2$ which we consider a given function. The distance
\begin{align}
    D(\partial \Omega^{\balpha}_h, \partial \Omega^{\balpha}) \equiv \left \| \mathbf{d}^{\balpha}_h - \mathbf{d}^{\balpha} \right \|_{L_2(\partial \hat{\Omega})}
\end{align}
serves as a measure of the approximation quality. \\
Likewise, we approximate the gradient by replacing $(u^{\balpha}, \mathbf{x}^{\balpha}) \rightarrow (u^{\balpha}_h, \mathbf{x}^{\balpha}_h)$ in (\ref{eq:grad_J_x}), i.e.,
\begin{align}
\label{grad_J_discrete}
    \frac{\mathrm{d} J}{\mathrm{d} \balpha} & \simeq \frac{\partial J_h^{\balpha}}{\partial u^{\balpha}_h} \left( \frac{\partial u^{\balpha}_h}{\partial \mathbf{x}^{\balpha}_h} \frac{\mathrm{d} \mathbf{x}^{\balpha}_h}{\mathrm{d} \balpha} + \frac{\partial u^{\balpha}_h}{\partial \balpha} \right) + \frac{\partial J_h^{\balpha}}{\partial \mathbf{x}^{\balpha}_h} \frac{\mathrm{d} \mathbf{x}^{\balpha}_h}{\mathrm{d} \balpha} + \frac{\partial J_h^{\balpha}}{\partial \balpha}.
\end{align}
At this point, it should be noted that for given $\balpha$, the exact evaluations of $J(\ldots)$ and the components of its gradient are independent of the particular choice of the coordinate system $\mathbf{x}^{\balpha}$. As such, the quality of the approximations introduced in (\ref{J_discrete}) and (\ref{grad_J_discrete}) depend solely on the numerical accuracy of $u^{\balpha}_h$, which in turn is affected by the parametric quality of $\mathbf{x}^{\balpha}_h$ and the distance of $\Omega_h^{\balpha}$ to the exact $\Omega^{\balpha}$. \\
We numerically treat (\ref{eq:Shape_Optimization}) based on a \textit{variable basis approach} (VBA) rather than a \textit{static basis approach} (SBA). In SBA, $u^{\balpha}_h$ and $\mathbf{x}^{\balpha}_h$ are constructed from the static tuple $(\mathcal{U}_h, \boldsymbol{\mathcal{V}_h})$, while in VBA the tuple $(\mathcal{U}_h^{\balpha}, \boldsymbol{\mathcal{V}}_h^{\balpha})$ may be chosen differently during each iteration and is tuned to the current needs. We make a choice based on the following principles:
\begin{enumerate}[label={\Alph*}., ref={\Alph*}]
    \item \label{enum:A} $\left \| \mathbf{d}_h^{\balpha} - \mathbf{d}^{\balpha} \right \|$, with $\mathbf{d}_h^{\balpha} \in \boldsymbol{\mathcal{V}}_h^{\balpha}$ is sufficiently small;
    \item \label{enum:B} $\mathbf{x}_h^{\balpha} \in \boldsymbol{\mathcal{V}}_h^{\balpha}$, resulting from $\mathbf{d}_h^{\balpha}$ in combination with (\ref{eq:x_weak_discretized}), is a bijection and preferably of high numerical quality;
    \item \label{enum:C} $u^{\balpha}_h \in \mathcal{U}_h^{\balpha}$ approximates $u^{\balpha}$ well.
\end{enumerate}
As such, for given $\balpha$, we select the tuple $(\mathcal{U}_h^{\balpha}, \boldsymbol{\mathcal{V}}_h^{\balpha})$ such that points \ref{enum:A} to \ref{enum:C} are satisfied with a minimal number of degrees of freedom (DOFs). \\
In SBA, a necessary condition for local optimality follows straightforwardly from the discretized counterpart of (\ref{eq:Shape_Optimization}) over the static tuple $(\mathcal{U}_h, \boldsymbol{\mathcal{V}}_h)$. In contrast, VBA necessitates basing such a condition on (\ref{eq:Shape_Optimization}) \textit{before} discretization. Hence, numerical assessment of local optimality in (\ref{eq:Shape_Optimization}) is obligatory, due to the approximate nature of $J_h^{\balpha}$ and its gradient. This may be regarded as a drawback since it can generate false positives / negatives caused by the truncation error at the current iterate. On the other hand, VBA allows for $\balpha$-specific feature-based basis selection for approximating both $\mathbf{x}^{\balpha}$ and $u^{\balpha}$, leading to a highly flexible scheme. When performing shape optimization in combination with EGG, a static $\mathcal{V}_h$ may be inappropriate for particular choices of $\balpha$ which results in grid-folding (impeding the evaluation of $J_h^{\balpha}$), hence justifying VBA-enabled feature-based basis selection in applications which are geometrically complex.
\begin{remark}
If we regard the truncation error $\tau(\balpha)$ in $u^{\balpha} = u^{\balpha}_h + \tau(\balpha)$ as a random variable drawn from some probability distribution, above methodology possesses many properties reminiscent of stochastic gradient descent \cite{bottou2010large}. As such, the convergence tolerance should be designed with the expected magnitude of $\tau(\balpha)$ (and its contribution to the gradient) in mind and hence taken generously. Here, we regard this as a minor shortcoming since we consider complex and highly nonconvex, nonlinear optimization problems in which the model error as well as the notorieties associated with nonconvex optimization (such as the danger of getting stuck in local minima) pose a greater threat to solution quality than the truncation error in practice. Furthermore, in most practical applications, a particular state vector need not be optimal in order to be considered adequate.
\end{remark}

\subsection{Gradient-Based Optimization Using an Adjoint Formulation}
In the following, we present a scheme that is suitable for gradient-based optimization, where all terms involved are assembled from expressions that have been derived fully symbolically. \\
For given $\balpha$, we assume that a suitable tuple $(\mathcal{U}_h^{\balpha}, \boldsymbol{\mathcal{V}}_h^{\balpha})$ has been chosen based on principles \ref{enum:A} to \ref{enum:C} (see section \ref{subsect:Discretization}). Particular methodologies for satisfying these principles depend on the application and are discussed in Section \ref{sect:Examples}. Reminiscent of (\ref{eq:test_space_x_f}), we introduce
\begin{align}
\label{eq:V_alphaf_h_U_alphaf_h}
    \boldsymbol{\mathcal{V}}^{\balpha,\mathbf{f}}_h \equiv \{ \mathbf{v} \in \boldsymbol{\mathcal{V}}^{\balpha}_h \enskip \vert \enskip \mathbf{v} = \mathbf{f} \text{ on } \partial \hat{\Omega} \} \quad \text{and} \quad \mathcal{U}^{\balpha,\mathbf{f}}_h \equiv \{ v \in \mathcal{U}^{\balpha}_h \enskip \vert \enskip v = f \text{ on } \partial \hat{\Omega} \},
\end{align}
where, as before, $\mathbf{f} \in \boldsymbol{\mathcal{V}}_h^{\balpha}$ and $f \in \mathcal{U}_h^{\balpha}$ by assumption. Equation (\ref{eq:V_alphaf_h_U_alphaf_h}) allows for the decomposition into \textit{boundary} ($\mathcal{B})$ and \textit{inner} ($\mathcal{I})$ bases:
\begin{align}
    \boldsymbol{\mathcal{V}}^{\balpha}_h = \boldsymbol{\mathcal{V}}^{\balpha,\mathcal{B}}_h \oplus \boldsymbol{\mathcal{V}}^{\balpha,\mathcal{I}}_h, \quad \text{with} \quad 
    \boldsymbol{\mathcal{V}}^{\balpha,\mathcal{I}}_h = \boldsymbol{\mathcal{V}}^{\balpha,\mathbf{0}}_h \quad \text{and} \quad     \boldsymbol{\mathcal{V}}^{\balpha,\mathcal{B}}_h = \boldsymbol{\mathcal{V}}^{\balpha}_h \setminus \boldsymbol{\mathcal{V}}^{\balpha,\mathcal{I}}_h,
\end{align}
and similarly for $\mathcal{U}_h^{\balpha} = \mathcal{U}^{\balpha,\mathcal{B}}_h \oplus \mathcal{U}^{\balpha,\mathcal{I}}_h$. Given $\mathbf{d}_h^{\balpha} \in \boldsymbol{\mathcal{V}}^{\balpha,\mathcal{B}}_h$ (see Section \ref{subsect:Discretization}), we introduce the mapping
\begin{align}
    \mathbf{x}^{\balpha}_h = \mathbf{x}_0^{\balpha} + \mathbf{d}_h^{\balpha}, \quad \text{with} \quad \mathbf{x}_0^{\balpha} = \sum_{\boldsymbol{\sigma}_i \in \boldsymbol{\mathcal{V}}^{\balpha,\mathcal{I}}_h} c_i^{\mathcal{I}} \boldsymbol{\sigma}_i \quad \text{and} \quad \mathbf{d}_h^{\balpha} = \sum_{\boldsymbol{\sigma}_i \in \boldsymbol{\mathcal{V}}^{\balpha,\mathcal{B}}_h} c_i^{\mathcal{B}} \boldsymbol{\sigma}_i,
\end{align}
where the $c^{\mathcal{B}}_i$ are known. We introduce the vector of weights $\mathbf{c}_{\mathcal{A}}$ (where the subscript $\mathcal{A}$ stands for \textit{'all'}), which is the concatenation of the vectors $\mathbf{c}_{\mathcal{I}}$ and $\mathbf{c}_{\mathcal{B}}(\balpha)$, containing the $c^{\mathcal{I}}_i$ and $c^{\mathcal{B}}_i$, respectively. Similarly, we introduce
\begin{align}
    u^{\balpha}_h = \sum_{\phi_i \in \mathcal{U}^{\balpha}_h} d_i \phi_i \quad \text{with the corresponding vector of weights} \quad (\mathbf{d}_{\mathcal{A}})_i = d_i.
\end{align}
With the introduction of the tuple $(\mathbf{c}_{\mathcal{A}}, \mathbf{d}_{\mathcal{A}})$, the discrete objective function is rewritten in the form
\begin{align}
    J_h^{\balpha}(\mathbf{x}_h^{\balpha}, u^{\balpha}_h, \balpha) \longrightarrow J_h^{\balpha} \left(\mathbf{c}_{\mathcal{A}}, \mathbf{d}_{\mathcal{A}}, \balpha \right), \quad \text{with} \quad \mathbf{c}_{\mathcal{A}} = \mathbf{c}_{\mathcal{A}}(\mathbf{c}_{\mathcal{I}}, \mathbf{c}_{\mathcal{B}}), \quad \mathbf{c}_{\mathcal{I}} = \mathbf{c}_{\mathcal{I}}(\mathbf{c}_{\mathcal{B}}) \quad \text{and} \quad \mathbf{d}_{\mathcal{A}} = \mathbf{d}_{\mathcal{A}}(\mathbf{c}_{\mathcal{A}}, \balpha).
\end{align}
With above concatenated dependencies in mind, the transposed gradient approximation reads:
\begin{align}
    \label{grad_Jh_c}
    \frac{\mathrm{d} J_h^{\balpha}}{\mathrm{d} \balpha}^T = \left[\frac{\mathrm{d} \cA}{\mathrm{d} \balpha} \right]^T \left( \left[\frac{\partial \dA}{\partial \cA} \right]^T \frac{\partial J_h^{\balpha}}{\partial \dA}^T + \frac{\partial J_h^{\balpha}}{\partial \cA}^T \right) + \left[ \frac{\mathrm{d} \dA}{\mathrm{d} \balpha} \right]^T \frac{\partial J_h^{\balpha}}{\partial \dA}^T + \frac{\partial J_h^{\balpha}}{\partial \balpha}^T,
\end{align}
where we denoted matrix quantities in square brackets. \\
Introducing the discrete EGG residual vector $\mathbf{F}^{\balpha}_h$ with
\begin{align}
    \left(\mathbf{F}^{\balpha}_h \right)_i = F(\mathbf{x}_h^{\balpha}, \boldsymbol{\sigma}_i), \quad \text{for} \quad \boldsymbol{\sigma}_i \in \boldsymbol{\mathcal{V}}^{\balpha,\mathcal{I}}_h \quad \text{and} \quad F(\cdot, \cdot) \text{ as defined in } (\ref{eq:F_operator}),
\end{align}
we use the implicit function theorem \cite{krantz2012implicit} to derive an expression for the gradient of $\cA$. We have:
\begin{align}
    \left[ \frac{\mathrm{d} \cA}{\mathrm{d} \balpha} \right]^T = \left[ \left[ \frac{\mathrm{d} \cI}{\mathrm{d} \balpha} \right]^T, \left[ \frac{\partial \cB}{\partial \balpha} \right]^T \right], \quad \text{with} \quad \left[ \frac{\mathrm{d} \cI}{\mathrm{d} \balpha} \right] = - \left[ \frac{ \partial \mathbf{F}^{\balpha}_h }{\partial \cI} \right]^{-1} \left[ \frac{ \partial \mathbf{F}^{\balpha}_h }{\partial \cB} \right] \left[ \frac{ \partial \cB }{\partial \balpha} \right].
\end{align}
Similarly, we define the residual vector $\mathbf{B}_h^{\balpha}$ of the discretized weak state equation (see (\ref{eq:u_weak_not_discretized})), with entries
\begin{align}
    \left( \mathbf{B}_h^{\balpha} \right)_i = B \left(u^{\balpha}_h, \mathbf{x}^{\balpha}_h, \balpha, \phi_i \right), \quad \text{for} \quad \phi_i \in \mathcal{U}^{\balpha, \mathcal{I}}_h.
\end{align}
The implicit function theorem yields
\begin{align}
    \left[ \frac{\partial \dA}{\partial \cA} \right] = - \left[ \frac{\partial \mathbf{B}_h^{\balpha}}{\partial \dA} \right]^{-1} \left[ \frac{\partial \mathbf{B}_h^{\balpha}}{\partial \cA} \right] \quad \text{and} \quad \left[ \frac{\mathrm{d} \dA}{\mathrm{d} \balpha} \right] = - \left[ \frac{\partial \mathbf{B}_h^{\balpha}}{\partial \dA} \right]^{-1} \left[ \frac{\partial \mathbf{B}_h^{\balpha}}{\partial \balpha} \right].
\end{align}
Substituting in (\ref{grad_Jh_c}) leads to
\begin{align}
\label{grad_Jh_ab}
    \frac{\mathrm{d} J_h^{\balpha}}{\mathrm{d} \balpha}^T = \left[\frac{\mathrm{d} \cA}{\mathrm{d} \balpha} \right]^T \mathbf{b} - \left[ \frac{\partial \mathbf{B}_h^{\balpha}}{\partial \balpha} \right]^T \mathbf{a} + \frac{\partial J_h^{\balpha}}{\partial \balpha}^T,
\end{align}
with
\begin{align}
\label{eq:a_b}
    \mathbf{a} = \left[ \frac{\partial \mathbf{B}_h^{\balpha}}{\partial \dA} \right]^{-T} \frac{\partial J_h^{\balpha}}{\partial \dA}^T \quad \text{and} \quad \mathbf{b} = - \left[ \frac{\partial \mathbf{B}_h^{\balpha}}{\partial \cA} \right]^T \mathbf{a} + \frac{\partial J_h^{\balpha}}{\partial \cA}^T.
\end{align}
Vector $\mathbf{a}$ is computed by solving the following linear system:
\begin{align}
\label{eq:a}
    \left[ \frac{\partial \mathbf{B}_h^{\balpha}}{\partial \dA} \right]^{T} \mathbf{a} = \frac{\partial J_h^{\balpha}}{\partial \dA}^T.
\end{align}
Vector $\mathbf{b}$ then follows from substituting $\mathbf{a}$ in (\ref{eq:a_b}). Furthermore, we have
\begin{align}
\label{eq:dcA_dbalpha_T_times_b_t}
    \left[\frac{\mathrm{d} \cA}{\mathrm{d} \balpha} \right]^T \mathbf{b} & = \left[ \left[ \frac{\mathrm{d} \cI}{\mathrm{d} \balpha} \right]^T, \left[ \frac{\partial \cB}{\partial \balpha} \right]^T \right] \begin{bmatrix} \mathbf{b}_{\mathcal{I}} \\[0.1cm] \mathbf{b}_{\mathcal{B}} \end{bmatrix} = \left[ \frac{ \partial \cB }{\partial \balpha} \right]^T \mathbf{q},
\end{align}
where
\begin{align}
    \mathbf{q} = - \left[ \frac{ \partial \mathbf{F}^{\balpha}_h }{\partial \cB} \right]^T \mathbf{e} + \mathbf{b}_{\mathcal{B}} \quad \text{with} \quad \mathbf{e} = \left[ \frac{ \partial \mathbf{F}^{\balpha}_h }{\partial \cI} \right]^{-T} \mathbf{b}_{\mathcal{I}}.
\end{align}
We compute the matrix-vector product
\begin{align}
\label{eq:compute_e}
    \mathbf{e} = \left[ \frac{ \partial \mathbf{F}^{\balpha}_h }{\partial \cI} \right]^{-T} \mathbf{b}_{\mathcal{I}} \quad \text{from the solution of} \quad \left[ \frac{ \partial \mathbf{F}^{\balpha}_h }{\partial \cI} \right]^{T} \mathbf{e} = \mathbf{b}_{\mathcal{I}}.
\end{align}
Finally, it should be noted that the matrix $\left[ \partial_{\balpha} \cB \right]$ in (\ref{eq:dcA_dbalpha_T_times_b_t}) depends on the collocation operator 
\begin{align}
    \pi^{\balpha}: \boldsymbol{\mathcal{V}} \setminus \boldsymbol{\mathcal{V}}^0 \rightarrow \boldsymbol{\mathcal{V}}^{\balpha, \mathcal{B}}_h, \quad \text{with} \quad \pi^{\balpha}(\mathbf{d}^{\balpha}) = \mathbf{d}_h^{\balpha}
\end{align}
and typically involves a sparse matrix - matrix inverse product. Upon transposing, the order of multiplication is reversed and the transposed inverse moves to the front. Therefore, we can treat matrix-vector products of the form $\left[ \partial_{\balpha} \cB \right]^{T} \mathbf{k}$ by inverting a sparse linear system and subsequent multiplication by a sparse matrix. We will present tangiable examples of this step in Section \ref{sect:Examples}. \\
In the following, we recapitulate all the necessary steps for computing the tuple $\left(J_h^{\balpha}, \mathrm{d}_{\balpha} J_h^{\balpha} \right)$ for given $\balpha$.
\begin{enumerate}[label=\textbf{S.\arabic*}:, ref=\textbf{S.\arabic*}]
    \item \label{enum:step_1} Choose an appropriate basis tuple $(\mathcal{U}_h^{\balpha}, \boldsymbol{\mathcal{V}}_h^{\balpha})$.
    \item \label{enum:step_2} Compute $\mathbf{d}_h^{\balpha}$ from $\mathbf{d}^{\balpha}$ using $\pi^{\balpha}$.
    \item \label{enum:step_3} Solve the nonlinear root-finding problem $\mathbf{F}_h^{\balpha}(\cI) = \mathbf{0}$, yielding the analysis-suitable mapping $\mathbf{x}^{\balpha}_h$.
    \item \label{enum:step_4} Solve the root-finding problem $\mathbf{B}_h^{\balpha} = \mathbf{0}$ using a suitable numerical algorithm. This yields the state variable $u^{\balpha}_h$.
    \item \label{enum:step_5} Substitute $(u^{\balpha}_h, \mathbf{x}^{\balpha}_h, \balpha)$ in $J(\cdot, \cdot, \cdot)$ to compute $J_h^{\balpha}$.
    \item \label{enum:step_6} Compute $\mathbf{a} \rightarrow \mathbf{b} \rightarrow \mathbf{e} \rightarrow \mathbf{q}$ and finally $\mathrm{d}_{\balpha} J_h^{\balpha}$ using (\ref{grad_Jh_ab}) to (\ref{eq:compute_e}).
\end{enumerate}
Due to the approximate nature of the tuple $(u^{\balpha}_h, \mathbf{x}^{\balpha}_h)$, we allow for a small amount of slack in the assessment of numerical feasibility, i.e., we replace
\begin{equation}
\begin{alignedat}{3}
    & g_i(u^{\balpha}, \mathbf{x}^{\balpha}, \balpha) \geq 0 \quad && \longrightarrow \quad && g_i(u^{\balpha}_h, \mathbf{x}^{\balpha}_h, \balpha) \geq \mu \quad \forall i \in \{1, \ldots, N_{\neq} \} \\
    & h_j(u^{\balpha}, \mathbf{x}^{\balpha}, \balpha) = 0 \quad && \longrightarrow \quad -\mu \leq \enskip && h_j(u^{\balpha}_h, \mathbf{x}^{\balpha}_h, \balpha) \leq \mu \quad \forall j \in \{1, \ldots, N_{=} \},
\end{alignedat}
\end{equation}
with $\mu > 0$ in (\ref{eq:Shape_Optimization}). The procedure that carries out $\ref{enum:step_1}$ to $\ref{enum:step_6}$, along with the relaxed constraints, is passed to a gradient-based optimization routine (such as IPOPT).

\subsection{Gradient Assembly Costs}
In the following, we analyse the computational costs of assembling the gradient. The majority of the costs result from assembling sparse matrices, as well as solving sparse linear systems, such as in (\ref{eq:a}). In order to compute $\mathbf{x}^{\balpha}_h$, we simultaneously assemble the quantities
\begin{align}
\label{eq:F_dF_tandem}
    \mathbf{F}^{\balpha}_h(\cI^i) \quad \text{and} \quad \left[ \frac{\partial \mathbf{F}^{\balpha}_h}{\partial \cI^i} \right]
\end{align}
during a joint element loop at the beginning of the $i$-th Newton iteration (see Section \ref{sect:EGG}). As such, this routine automatically yields the matrix $\left[ \partial_{\cI} \mathbf{F}^{\balpha}_h \right]$ at the last step. In the following, we assume that $\partial \Omega^{\balpha}_D = \emptyset$ for convenience, i.e., $u^{\balpha}_h$ is free of Dirichlet data or the data is enforced with Nitsche's method \cite{hansbo2002unfitted}. If $\mathbf{B}_h^{\balpha}$ is linear, assembling $\left[ \partial_{\dA}\mathbf{B}^{\balpha}_h \right]$ is a precursor to computing $u^{\balpha}_h$ and hence available. If $\mathbf{B}_h^{\balpha}$ is nonlinear, we recommend basing an iterative algorithm on Newton's method and computing the residual and its derivative in tandem, as in (\ref{eq:F_dF_tandem}). As such, additional cost factors are assembling $\left[ \partial_{\cA} \mathbf{B}_h^{\balpha} \right]$ and solving a number of sparse linear equations. Due to the nonlinear nature of $\mathbf{F}^{\balpha}_h$, the cost of computing $\mathrm{d}_{\balpha} J^{\balpha}_h$ is of the same order as a discrete evaluation of $J(\cdot, \cdot, \cdot)$ (regardless of the length of $\balpha$). \\
Finally, we note that the discrete constraint gradients associated with the $g_i$ and $h_j$ are efficiently computed by replacing $J_h^{\balpha}$ by the corresponding term in (\ref{grad_Jh_c}) and repeating steps (\ref{grad_Jh_ab}) to (\ref{eq:compute_e}). Hereby the required matrices can be reused from the assembly of the gradient.

\subsection{Memory-Saving Strategies in Large-Scale Applications}
In light of enabling large-scale optimization as well as the prospect of extending the presented methodology to volumetric applications, in the following, we discuss ways to avoid the memory-consuming assembly of the matrices involved in computing $u^{\balpha}_h$, $\mathbf{x}^{\balpha}_h$ and $J^{\balpha}_h$. \\
Memory-saving strategies are based on the observation that matrices only appear in the form of matrix-vector products during the assembly of $\mathrm{d}_{\balpha} J^{\balpha}_h$. Let $\mathbf{B} = \mathbf{B}(\ldots, \mathbf{q}, \ldots)$. Then, we have
\begin{align}
\label{eq:Gateaux_approx}
    \left[ \frac{\partial \mathbf{B}(\ldots, \mathbf{q}, \ldots)}{\partial \mathbf{q}} \right] \mathbf{a} \simeq \frac{ B(\ldots, \mathbf{q} + \epsilon \mathbf{a}, \ldots) - B(\ldots, \mathbf{q}, \ldots) }{\epsilon},
\end{align}
for $\epsilon > 0$ small. As such, in steps (\ref{grad_Jh_ab}) to (\ref{eq:compute_e}), matrix-vector products can be approximated using (\ref{eq:Gateaux_approx}). Since Krylov-subspace (KS) methods such as GMRES \cite{saad1986gmres} only require matrix-vector products, we combine a KS-method with (\ref{eq:Gateaux_approx}) for solving linear systems as they appear in, e.g., equation (\ref{eq:a}). Reminiscent of \textit{Newton-Krylov} \cite{knoll2004jacobian}, this principle may be extended to the computation of $u^{\balpha}_h$ and $\mathbf{x}^{\balpha}_h$, hence completely bypassing matrix assembly in steps \ref{enum:step_1} to \ref{enum:step_6}. Hereby, we regard the cumulative error contribution to $\mathrm{d}_{\balpha} J^{\balpha}_h \simeq \mathrm{d}_{\balpha} J^{\balpha}$ as negligible compared to other sources (such as the model error and the truncation resulting from the numerical scheme). The optimal choice of $\epsilon$ is discussed in \cite{knoll2004jacobian}.
\section{Examples}
\label{sect:Examples}
In this section we apply the methodology from Section \ref{subsect:Discretization} to selected test cases. We consider the first example a validation test case, in which the exact minimizer can be computed exactly (up to machine precision). Hereby, we compare the results of VBA (see Section \ref{sect:Computational_Approach}) to the exact minimum. Furthermore, we compare the VBA results to those resulting from taking the tuple $(\mathcal{U}^{\balpha}_h, \boldsymbol{\mathcal{V}}^{\balpha}_h)$ static (SBA). In the second case, we consider the design of a cooling element, whereby the plausibility of the outcome can only be assessed using physical reasoning. \\
Both examples have been carefully selected in order to be geometrically challenging. We implemented the scheme from Section \ref{sect:Computational_Approach} in the open-source Python library \textit{Nutils} \cite{gertjan_van_zwieten_2019_3243447}.

\subsection{A Validation Example with Known Exact Solution}
\label{sect:Validation}
\begin{figure}[h!]
  \centering
  \includegraphics[width=.9\linewidth]{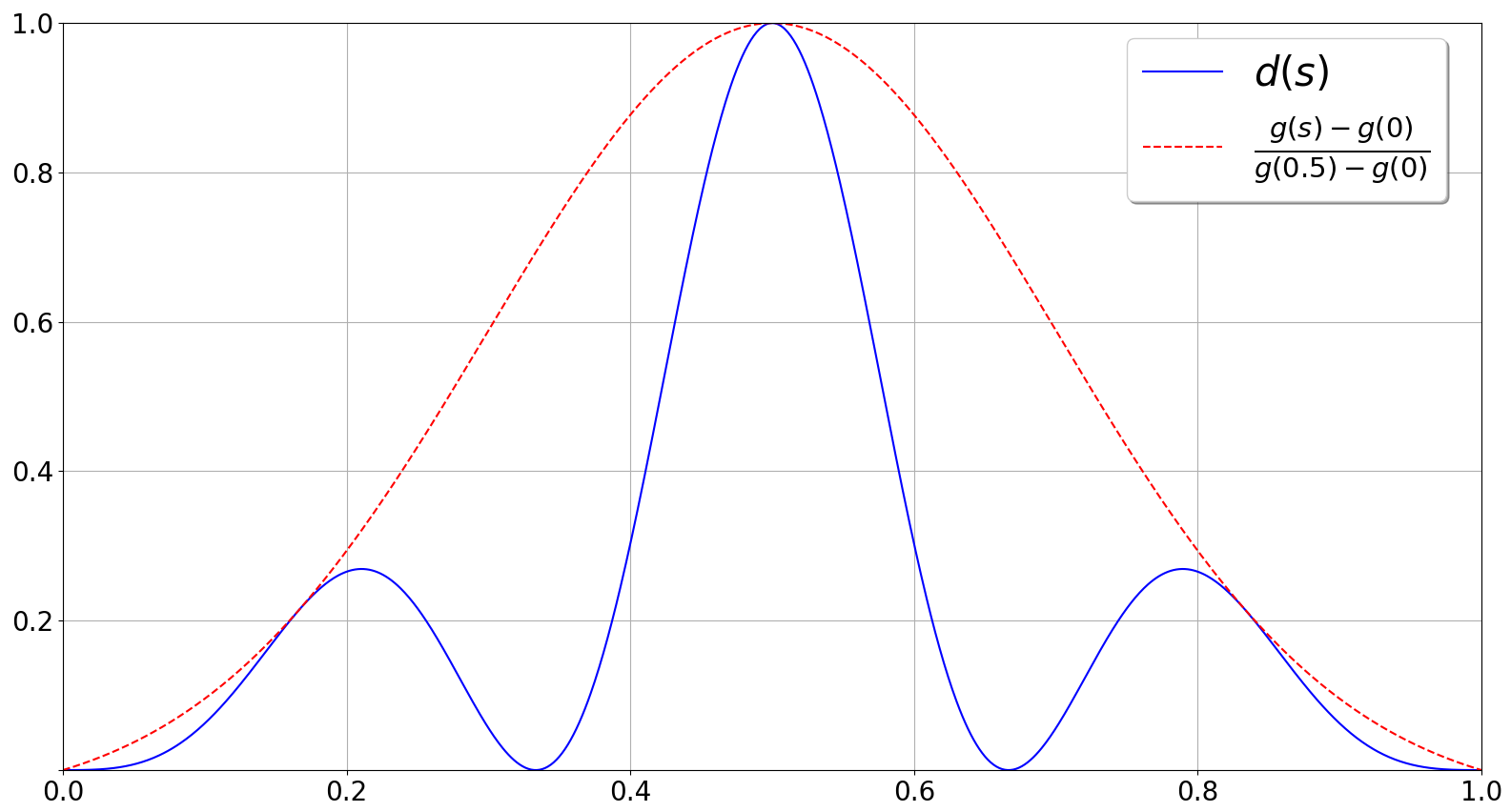}
  \caption{Plot of $d(s)$ (see (\ref{eq:validation_d})) and the corresponding envelope function over the interval $s \in [0, 1]$.}
  \label{fig:validation_d_envelope}
\end{figure}

\begin{figure}[h!]
  \centering
  \includegraphics[width=.5\linewidth]{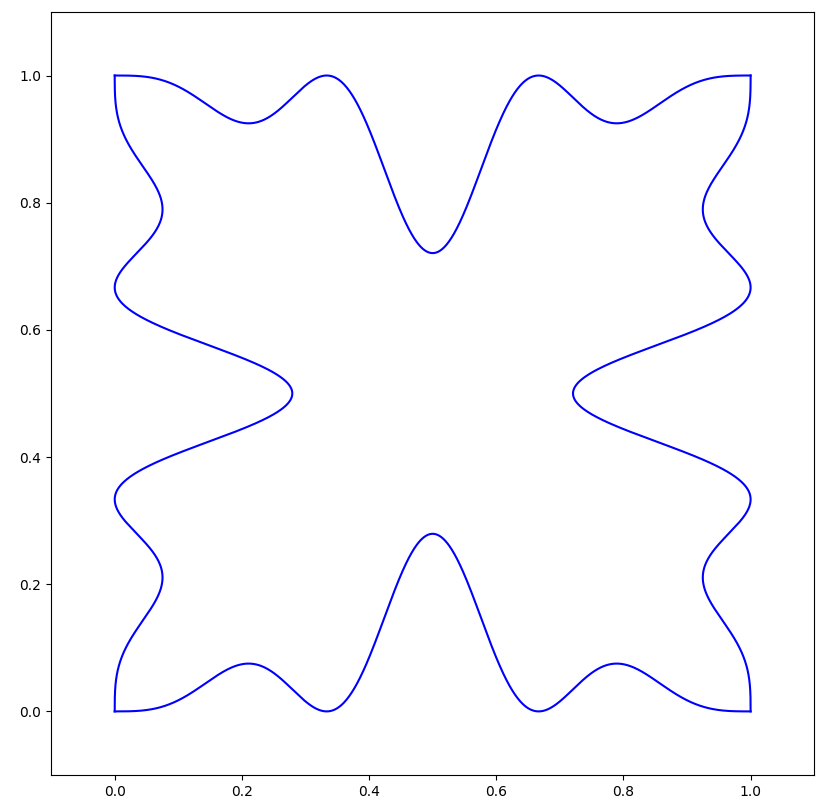}
  \caption{The contours of the domain $\Omega^{\balpha*}$ that corresponds to the exact minimizer $\balpha*$.}
  \label{fig:validation_d_exact}
\end{figure}

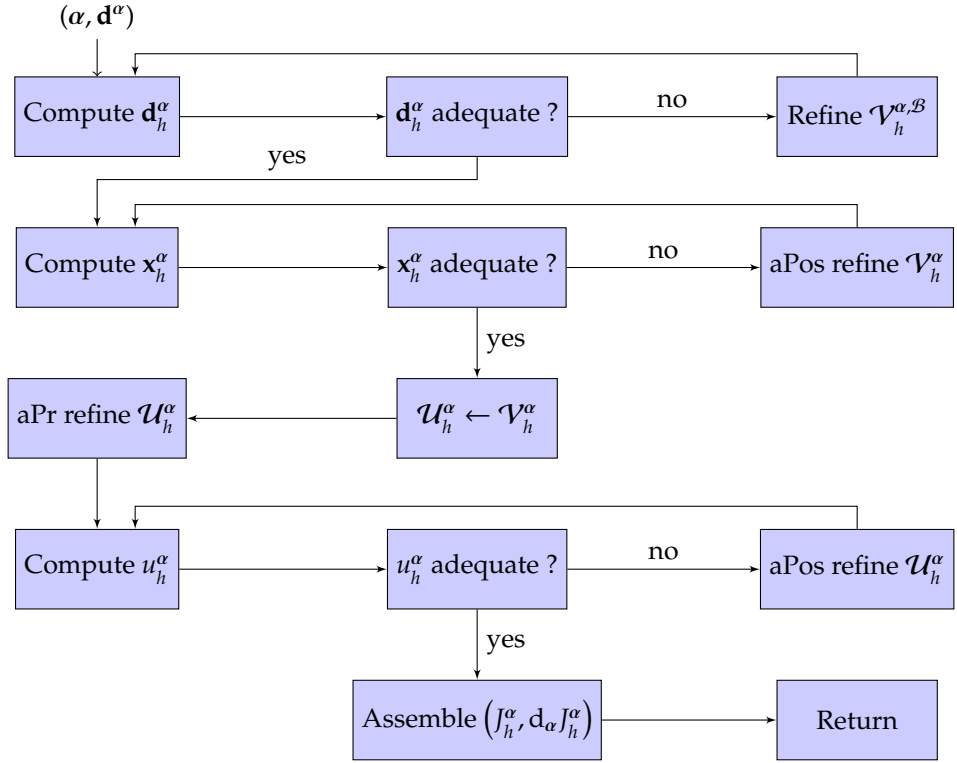
\begin{figure}[h!]
\centering
\begin{tikzpicture}[auto, node distance=2cm,>=latex']
    \node [block, pin={[pinstyle]above:$(\balpha, \mathbf{d}^{\balpha})$}] (step1) {Compute $\mathbf{d}^{\balpha}_h$};
    \node [block, right of=step1, node distance=5cm] (step1_adequate) {$\mathbf{d}^{\balpha}_h$ adequate ?};
    \node [block, right of=step1_adequate, node distance=5cm] (refine_dB) {Refine $\mathcal{V}^{\balpha, \mathcal{B}}_h$};
    \node [block, below of=step1, node distance=2cm] (step2) {Compute $\mathbf{x}^{\balpha}_h$};
    \node [block, right of=step2, node distance=5cm] (step2_adequate) {$\mathbf{x}^{\balpha}_h$ adequate ?};
    \node [block, right of=step2_adequate, node distance=5cm] (refine_V) {aPos refine $\mathcal{V}^{\balpha}_h$};
    \node [block, below of=step2_adequate, node distance=2cm] (step3) {$\mathcal{U}^{\balpha}_h \leftarrow \mathcal{V}^{\balpha}_h$};
    \node [block, left of=step3, node distance=5cm] (aPr_U) {aPr refine $\mathcal{U}^{\balpha}_h$};
    \node [block, below of=aPr_U, node distance=2cm] (step4) {Compute $u^{\balpha}_h$};
    \node [block, right of=step4, node distance=5cm] (step4_adequate) {$u^{\balpha}_h$ adequate ?};
    \node [block, right of=step4_adequate, node distance=5cm] (refine_U) {aPos refine $\mathcal{U}^{\balpha}_h$};
    \node [block, below of=step4_adequate, node distance=2cm] (compute_J) {Assemble $\left(J^{\balpha}_h, \mathrm{d}_{\balpha} J^{\balpha}_h \right)$};
    \node [block, right of=compute_J, node distance=5cm] (return) {Return};

    \draw [->] (step1) -- node {} (step1_adequate);
    \draw [->] (step1_adequate) -- node {no} (refine_dB);
    \path[draw,->] 
    (refine_dB.north) -- ++(0,0.3cm) -- ++(-9.5cm,0) -- ++(0,-0.3cm);
    \path[draw,->] 
    (step1_adequate.south) -- ++(0,-0.3cm) -| node[pos=0.25, above] {yes} (step2.north);
    \draw [->] (step2) -- node {} (step2_adequate);
    \draw [->] (step2_adequate) -- node {no} (refine_V);
    \path[draw,->] 
    (refine_V.north) -- ++(0,0.3cm) -- ++(-9.5cm,0) -- ++(0,-0.3cm);
    \path[draw,->] 
    (step2_adequate.south) -- node {yes} (step3.north);
    \draw [->] (step3) -- node {} (aPr_U);
    \draw [->] (aPr_U) -- node {} (step4);
    \draw [->] (step4) -- node {} (step4_adequate);
    \draw [->] (step4_adequate) -- node {no} (refine_U);
    \path[draw,->] 
    (refine_U.north) -- ++(0,0.3cm) -- ++(-9.5cm,0) -- ++(0,-0.3cm);
    \draw [->] (step4_adequate) -- node {yes} (compute_J);
    \draw [->] (compute_J) -- node {} (return);
\end{tikzpicture}
\caption{Block diagram summarizing all the steps required for computing the tuple $\left(J^{\balpha}_h, \mathrm{d}_{\balpha} J^{\balpha}_h \right)$.}
\label{fig:block_diagram_J}
\end{figure}

We are considering the example of a domain fenced-off by four parametric curves that are given by an envelope function multiplied by a cosine, whereby the amplitude of the cosine is a degree of freedom in $\balpha = (\alpha_1, \alpha_2, \alpha_3, \alpha_4)$. We define the function
\begin{align}
\label{eq:validation_d}
    d(s) = \underbrace{\left( \frac{g(s) - g(0)}{g(0.5) - g(0)} \right)}_{\text{envelope function}} \underbrace{\left(\frac{1 - \cos(\omega \pi s)}{2} \right)}_{\text{trigonometric component}}, \quad \text{where} \quad g(s) = \mathrm{exp}\left(- \frac{\left(s - \tfrac{1}{2} \right)^2}{2 \sigma^2} \right),
\end{align}
with $(\omega, \sigma) = (6, 0.2)$. Note that $d(0) = d(1) = 0$ and $d(0.5) = 1$, $d(s)$ and the envelope function are depicted in Figure \ref{fig:validation_d_envelope}.
\noindent Let $\partial \hat{\Omega} = \bar{\gamma}_{S} \cup \bar{\gamma}_{E} \cup \bar{\gamma}_{N} \cup \bar{\gamma}_{W}$, where the $\gamma_{\beta}, \enskip \beta \in \{S, E, N, W\}$ denote the southern, eastern, northern, and western boundary of $\partial \hat{\Omega}$, respectively. The contour function $\mathbf{d}^{\balpha}: \partial \hat{\Omega} \rightarrow \partial \Omega^{\balpha}$ reads:
\begin{align}
    \mathbf{d}^{\balpha}(\boldsymbol{\xi}) = \left\{
        \begin{array}{ll}
            \left( \xi_1, \alpha_1 d(\xi_1) \right)^T & \quad \boldsymbol{\xi} \in \bar{\gamma}_{S} \\
            \left( 1 - \alpha_2 d(\xi_2), \xi_2 \right)^T & \quad \boldsymbol{\xi} \in \bar{\gamma}_{E} \\
            \left( \xi_1, 1 - \alpha_3 d(\xi_1) \right)^T & \quad \boldsymbol{\xi} \in \bar{\gamma}_{N} \\
            \left( \xi_2, \alpha_4 d(\xi_2) \right)^T & \quad \boldsymbol{\xi} \in \bar{\gamma}_{W}
        \end{array}
    \right.,
\end{align}
while
\begin{align}
    \boldsymbol{\mathcal{\lambda}} = \{ \boldsymbol{\alpha} \in \mathbb{R}^4 \enskip \vert \enskip \mathbf{0} \leq \boldsymbol{\alpha} \leq \tfrac{2}{5} \mymathbb{1} \},
\end{align}
where $\mymathbb{1}$ is a vector of ones. \\
Here, we base $\pi^{\balpha}: \boldsymbol{\mathcal{V}} \setminus \boldsymbol{\mathcal{V}}^{\mathbf{0}} \rightarrow \boldsymbol{\mathcal{V}}^{\balpha, \mathcal{B}}_h$ (see Section \ref{sect:Computational_Approach}) on an $L_2(\partial \hat{\Omega})$ projection. For given $\boldsymbol{\mathcal{V}}^{\balpha, \mathcal{B}}_h$, we hence have
\begin{align}
\label{eq:L2_ddomain}
    \left[ \frac{\partial \cB}{\partial \balpha} \right]^T = - \left[ \frac{\partial \mathbf{D}^{\balpha}}{\partial \balpha} \right]^T \left[ \frac{\partial \mathbf{D}^{\balpha}}{\partial \cB}  \right]^{-T}, \quad \text{where} \quad \mathbf{D}_i^{\balpha} = \int_{\partial \hat{\Omega}} \boldsymbol{\sigma}_i^{\mathcal{B}} \cdot \left( \mathbf{d}^{\balpha} - \mathbf{d}_h^{\balpha}(\cB) \right) \mathrm{d}\gamma, \quad \text{with} \quad \boldsymbol{\sigma}_i^{\mathcal{B}} \in \boldsymbol{\mathcal{V}}^{\balpha, \mathcal{B}}_h.
\end{align}
In (\ref{eq:L2_ddomain}), we take matrix-vector products in the same way as in Section \ref{sect:Computational_Approach}. \\
We base our state variable residual on the following PDE-problem:
\begin{align}
\label{eq:validation_u_eq}
    -\Delta u^{\balpha} = -\Delta f^{\balpha}, \quad \text{s.t.} \quad u^{\balpha}\vert_{\partial \hat{\Omega}} = f^{\balpha} \quad \text{where} \quad f^{\balpha} = \det [J^{\balpha}] = \det \left[ \frac{\partial \mathbf{x}^{\balpha}}{\partial \boldsymbol{\xi}} \right].
\end{align}
Clearly, the exact solution of (\ref{eq:validation_u_eq}) satisfies $u^{\balpha} = \det [J^{\balpha}]$. We derive the weak form of (\ref{eq:validation_u_eq}) and implement the boundary conditions using Nitsche's method. This leads to
\begin{align}
    \left( \mathbf{B}^{\balpha}_h \right)_i = \left( \nabla(u^{\balpha}_h - f^{\balpha}), \nabla \phi_i \right)_{\Omega^{\balpha}_h} - \int \limits_{\partial \Omega^{\balpha}_h} \phi_i \frac{\partial u^{\balpha}_h}{\partial \mathbf{n}} \mathrm{d} \gamma - \int \limits_{\partial \Omega^{\balpha}_h} (u^{\balpha}_h - f^{\balpha}) \frac{\partial \phi_i}{\partial \mathbf{n}} \mathrm{d} \gamma + \eta_i \int \limits_{\partial \Omega^{\balpha}_h} (u^{\balpha}_h - f^{\balpha}) \phi_i \mathrm{d} \gamma,
\end{align}
with $\phi_i \in \mathcal{U}^{\balpha}_h$. Here, $\partial / \partial \mathbf{n}$ denotes the outward normal derivative with respect to $\Omega^{\balpha}_h$ and $\eta_i \gg 1$ is a penalty parameter. We use
\begin{align}
    \eta_i = \left\{
        \begin{array}{ll}
            c I_i^{-1} & \quad I_i > 0 \\
            0 & \quad \text{else}
        \end{array}
    \right., \quad \text{where} \quad I_i = \int_{\partial \Omega^{\balpha}_h} \phi_i \mathrm{d} \gamma \quad \text{and} \quad c = 10^3.
\end{align}
The objective function reads:
\begin{align}
    J^{\balpha} = \left \| u^{\balpha} \right \|_{\hat{\Omega}} + \frac{1}{2} \left \| \balpha \right \|^2 \implies J^{\balpha}_h = \left \| u^{\balpha}_h \right \|_{\hat{\Omega}} + \frac{1}{2} \left \| \balpha \right \|^2
\end{align}
and there are no further constraints. Since $u^{\balpha} = \det [J^{\balpha}]$, we have
\begin{align}
    \left \| u^{\balpha} \right \|_{\hat{\Omega}} = \operatorname{Area}(\Omega^{\balpha}) = 1 + \sum_i \alpha_i A, \quad \text{where} \quad A = -\int_{[0, 1]} d(s) \mathrm{d}s.
\end{align}
We compute the exact value of $A$ up to machine precision, which yields $A \simeq 0.2374$. The exact minimum over $\balpha \in \boldsymbol{\mathcal{\lambda}}$ is assumed at $\balpha* = A {\mymathbb{1}}$ and yields $J^{\balpha *} \simeq 0.8873$. The contours of the resulting domain $\Omega^{\balpha *}$ are depicted in Figure \ref{fig:validation_d_exact}.
 \\
For increasingly fine $(\mathcal{U}^{\balpha}_h, \boldsymbol{\mathcal{V}}^{\balpha}_h)$, the minimum of the discretized optimization problem should converge to the exact minimum, allowing us to test the consistency of the scheme. \\
In the following, we discuss how to choose the tuple $(\mathcal{U}^{\balpha}_h, \boldsymbol{\mathcal{V}}^{\balpha}_h)$ during each iteration. We start by dividing $\hat{\Omega}$ into a structured set of elements, resulting from the bivariate knot vector $\boldsymbol{\Xi}^{p_1, p_2} = \Xi^{p_1} \times \mathcal{H}^{p_2}$, where the $p_i$ denote the order. Here, we restrict ourselves to bicubic bases, i.e., $p_1 = p_2 = 3$. With points \ref{enum:A} to \ref{enum:C} (see Section \ref{sect:Computational_Approach}) in mind, we repeatedly refine the $\phi_i \in \mathcal{V}^{\balpha, \mathcal{B}}_h \subset \mathcal{V}^{\balpha}_h$, where $\mathcal{V}^{\balpha}_h$ is initialized to the coarse-grid basis resulting from $\boldsymbol{\Xi}^{p_1, p_2}$. Let the $i$-th contribution to the projection residual be denoted by $r_i(\mathbf{d}^{\balpha}_h)$, where
\begin{align}
\label{eq:projection_residual}
    R(\mathbf{d}^{\balpha}_h)^2 = \frac{1}{2} \int \limits_{\partial \hat{\Omega}} \left \| \mathbf{d}^{\balpha} - \mathbf{d}^{\balpha}_h \right \|^2 \mathrm{d} \gamma = \frac{1}{2} \sum \limits_{\phi_i \in \mathcal{V}^{\balpha}_h} \int \limits_{\partial \hat{\Omega}} \phi_i \left \| \mathbf{d}^{\balpha} - \mathbf{d}^{\balpha}_h \right \|^2 \mathrm{d} \gamma \equiv \frac{1}{2} \sum_i r_i^2(\mathbf{d}^{\balpha}_h).
\end{align}
We refine $\phi_i \in \mathcal{V}^{\balpha, \mathcal{B}}_h$ whenever $r_i(\mathbf{d}^{\balpha}_h)$ exceeds a threshold $\mu_i$. The threshold is of the form
\begin{align}
\label{eq:mu_i}
    \mu_i = \frac{\mu}{ \sqrt{\left\| \phi_i \right\|_{L_2(\partial \hat{\Omega})}}},
\end{align}
where $\mu$ is a small positive constant that tunes the accuracy of $\mathbf{d}^{\balpha}_h$. Note that in (\ref{eq:projection_residual}), we have made use of the partition of unity property that holds for THB-spline bases. \\
As a next step, we compute $\mathbf{x}^{\balpha}_h$ using the methodology from Section \ref{sect:EGG}. As the choice of $\boldsymbol{\mathcal{V}}^{\balpha}_h$, at this point, is solely based on accurately resolving the boundary contours, it may be too optimistic (in terms of the number of inner DOFs) for computing a folding-free mapping. In the case of folding, we apply a posteriori refinement to \textit{defective} elements (i.e., elements $\epsilon \subset \hat{\Omega}$ on which $\det [J^{\balpha}_h](\boldsymbol{\xi}_i) < 0$ for some $\boldsymbol{\xi}_i \in \epsilon$) by refining all $\phi \in \mathcal{V}^{\balpha}_h$ that are non-vanishing on that element. The defective mapping is prolonged to the refined space and serves as an initial guess for recomputing it from the enriched space. This step may be repeated until $\boldsymbol{\mathcal{V}}^{\balpha}_h$ is such that $\mathbf{x}^{\balpha}_h \in \boldsymbol{\mathcal{V}}^{\balpha}_h$ is folding-free. Note that, although the proposed methodology is robust in practice, it may lead to over-refinement. The methodology may be combined with the refinement strategies proposed in \cite{hinz2020goaloriented}, which avoid over-refinement.
\begin{remark}
Here, we base the selection of $\boldsymbol{\mathcal{V}}^{\balpha}_h$ on a posteriori strategies, which necessitates recomputing $\mathbf{x}^{\balpha}_h$ after each refinement. Choosing the coarse grid basis properly (i.e., not too coarse), we typically did not encounter more than $1-2$ a posteriori refinements in the cases considered in this work. Fortunately, the defective mappings can be used as an initial guess for the recomputed one, significantly reducing computational costs. \\
Reliable a priori refinement strategies are however desirable and constitute a topic for future research.
\end{remark}
After achieving bijectivity, additional refinement can be applied in order to further-improve the quality of the mapping. A posteriori strategies that rely on the \textit{Winslow} functional \cite{charakhch1997variational} are discussed in \cite{hinz2020goaloriented}. \\
Upon completion, we are in the possession of an analysis-suitable $\mathbf{x}^{\balpha}_h: \hat{\Omega} \rightarrow \Omega^{\balpha}_h$ from the appropriately refined $\boldsymbol{\mathcal{V}}^{\balpha}_h$. As a next step, we choose a suitable space $\mathcal{U}^{\balpha}_h$. Heuristically, there exists a strong correlation between the regions in which $\mathcal{V}^{\balpha}_h$ has been refined in order to yield an analysis-suitable $\mathbf{x}^{\balpha}_h$ and the regions that ought to be refined in order to accurately approximate $u^{\balpha}$. As such, we initialize $\mathcal{U}^{\balpha}_h$ to the current choice of $\mathcal{V}^{\balpha}_h$. In our current implementation we always base $\mathcal{U}^{\balpha}_h$ on $\mathcal{V}^{\balpha}_h$ or a (possibly repeated) uniform h-refinement thereof. However, for more flexibility, we briefly recapitulate possible feature-based refinement strategies. \\
In all cases, plausible a priori (aPr) strategies refine elements that are too large (on $\Omega^{\balpha}_h$), while a posterori (aPos) strategies depend on the underlying PDE-problem. In the case of (\ref{eq:validation_u_eq}), aPos-refinement can be based on the \textit{strong residual norm}
\begin{align}
\label{eq:strong_residual_norm}
    m_{\operatorname{SR}} = \int \limits_{\Omega^{\balpha}_h} \left( \Delta u^{\balpha}_h - \Delta \det [J^{\balpha}_h] \right)^2 \mathrm{d}S + \mathcal{F}(u^{\balpha}_h \vert_{\partial \hat{\Omega}}),
\end{align}
where $\mathcal{F}(\cdot): \mathcal{U}^{\balpha, \mathcal{B}}_h \rightarrow \mathbb{R}^{+}$ is a suitably-chosen penalty term that gauges how well the boundary condition is resolved by Nitsche's method. Note that (\ref{eq:strong_residual_norm}) requires taking $\mathbf{x}^{\balpha}_h \in C^2(\hat{\Omega})$, which is satisfied if we utilize bicubics with maximum regularity. Equation (\ref{eq:strong_residual_norm}) may then be decomposed into basis function wise contributions (as in (\ref{eq:projection_residual})) or serve as a cost function for \textit{dual weighted residual} (DWR) based aPos refinement \cite{rannacher2004adaptive}. \\
Alternatively, the \textit{weak residual norm} $m_{\operatorname{WR}}$, with
\begin{align}
    m_{\operatorname{WR}} = \sum_i c_{\operatorname{WR, i}}^2 \quad \text{and} \quad c_{\operatorname{WR, i}} = B(u^{\balpha}_h, \mathbf{x}^{\alpha}_h, \balpha, \psi_i)
\end{align}
may be utilized. Here, the $\psi_i$ are taken from a space $\bar{\mathcal{U}} \supset \mathcal{U}^{\balpha}_h$ that results from uniformly refining $\mathcal{U}^{\balpha}_h$ in $p$ or $h$. For more details, we refer to \cite{gravesen2012planar}. \\
Upon completion of an adequate state variable approximation $u^{\balpha}_h$, we are in the position to assemble the tuple $\left(J^{\balpha}_h, \mathrm{d}_{\balpha} J^{\balpha}_h \right)$ utilizing the principles from Section \ref{sect:Computational_Approach}. All the required steps are summarized in Figure \ref{fig:block_diagram_J}.
 \\
In the following, we present the results of a computational approach for various values of $\mu$ (see equation (\ref{eq:mu_i})) and $u_{\text{ref}}$, where $u_{\text{ref}}$ refers to the number of aPr $h$-refinements of $\mathcal{U}^{\balpha}_h$ with respect to $\mathcal{V}^{\balpha}_h$. For this, the procedure that corresponds to Figure \ref{fig:block_diagram_J} has been passed to a SLSQP \cite{kraft1988software} routine. In all cases, we use the initial guess $\balpha_0 = \mathbf{0}$.
\begin{table}[h!]
\centering
\subfloat[Table showing $\vert \min J^{\balpha}_h - \min J^{\balpha} \vert$.]{
\begin{tabular}{l|c|c|c}
  \backslashbox{\small $\mu$}{\small $u_{\text{ref}}$}  & $0$ & $1$ & $2$ \\ \hline
$10^{-2}$ & $8.2 \times 10^{-3}$ & $4.6 \times 10^{-4}$ & $7.7 \times 10^{-6}$  \\ 
$10^{-4}$ & $3.7 \times 10^{-3}$ & $1.6 \times 10^{-4}$ & $5.5 \times 10^{-6}$  \\
$10^{-6}$ & $2.7 \times 10^{-4}$ & $1.3 \times 10^{-5}$ & $5.2 \times 10^{-7}$
\end{tabular} } \qquad
\subfloat[$\#$iterations required until convergence.]{
\begin{tabular}{l|c|c|c}
  \backslashbox{\small $\mu$}{\small $u_{\text{ref}}$}  & $0$ & $1$ & $2$ \\ \hline
$10^{-2}$ & $4$ & $2$ & $3$  \\ 
$10^{-4}$ & $4$ & $3$ & $3$  \\
$10^{-6}$ & $2$ & $3$ & $3$
\end{tabular} }
\caption{Tables showing $\vert \min J^{\balpha}_h - \min J^{\balpha} \vert$ (a) and the required number of iterations until convergence is reached (b) for various combinations of $(\mu, u_{\text{ref}})$.}
\label{tab:minJ_mu_uref}
\end{table}

\begin{table}[h!]
\centering
\subfloat[Average $\#$DOFs for $\mathbf{x}^{\balpha}_h$.]{
\begin{tabular}{l|c|c|c}
  \backslashbox{\small $\mu$}{\small $u_{\text{ref}}$}  & $0$ & $1$ & $2$ \\ \hline
$10^{-2}$ & $483.5$ & $389$ & $452$  \\ 
$10^{-4}$ & $735.5$ & $676$ & $676$  \\
$10^{-6}$ & $1145$ & $1460$ & $1460$
\end{tabular}} \qquad
\subfloat[Average $\#$DOFs for $u^{\balpha}_h$.]{
\begin{tabular}{l|c|c|c}
  \backslashbox{\small $\mu$}{\small $u_{\text{ref}}$}  & $0$ & $1$ & $2$ \\ \hline
$10^{-2}$ & $241.75$ & $625$ & $2641$  \\ 
$10^{-4}$ & $367.75$ & $1169$ & $4337$  \\
$10^{-6}$ & $572.5$ & $2737$ & $10609$
\end{tabular} }
\caption{Tables showing the average of the number of DOFs involved in computing $\mathbf{x}^{\balpha}_h$ (a) and $u^{\balpha}_h$ (b) (over all iterations) for various combinations of $(\mu, u_{\text{ref}})$.}
\label{tab:ndofs_mu_uref}
\end{table}
Table \ref{tab:minJ_mu_uref} shows the discrepancy between the exact objective function minimum and its numerical approximation $|\min J^{\balpha}_h - \min J^{\balpha} \vert$ (a) and the required number of iterations until convergence is reached (b), for all possible combinations of $\mu = (10^{-2}, 10^{-4}, 10^{-6})$ and $u_{\text{ref}} = (0, 1, 2)$. In all cases, we initialized $\mathcal{V}^{\balpha}_h$ to a bicubic B-spline basis with $7$ elements per coordinate direction and maximum regularity.
\begin{table}[h!]
\centering
\subfloat[Table showing $\vert \min J^{\balpha}_h - \min J^{\balpha} \vert$. \label{tab:minJ_nelems_uref_min}]{
\begin{tabular}{l|c|c|c}
  \backslashbox{\small $n_e$}{\small $u_{\text{ref}}$}  & $0$ & $1$ & $2$ \\ \hline
$12$ & $8.4 \times 10^{-3}$ & $5.4 \times 10^{-4}$ & $6.5 \times 10^{-6}$  \\ 
$16$ & $7.2 \times 10^{-3}$ & $4.1 \times 10^{-4}$ & $1.1 \times 10^{-5}$  \\
$23$ & $4.7 \times 10^{-3}$ & $2.4 \times 10^{-4}$ & $8.3 \times 10^{-6}$
\end{tabular} } \qquad
\subfloat[$\#$iterations required until convergence. \label{tab:minJ_nelems_uref_niter}]{
\begin{tabular}{l|c|c|c}
  \backslashbox{\small $n_e$}{\small $u_{\text{ref}}$}  & $0$ & $1$ & $2$ \\ \hline
$12$ & $4$ & $3$ & $2$  \\ 
$16$ & $4$ & $3$ & $2$  \\
$23$ & $4$ & $3$ & $2$
\end{tabular} }
\caption{Tables showing $\vert \min J^{\balpha}_h - \min J^{\balpha} \vert$ (a) and the required number of iterations until convergence is reached (b) in the SBA case for various combinations of $(n_e, u_{\text{ref}})$.}
\label{tab:minJ_nelems_uref}
\end{table}
Finally, Tables \ref{tab:minJ_nelems_uref} and \ref{tab:ndofs_nelems_uref} show the corresponding results from a static basis approach. Hereby $n_e$ denotes the number of elements we used per coordinate direction. They have been carefully selected to yield roughly the same number of DOFs for both $\mathbf{x}^{\balpha}_h$ and $u^{\balpha}_h$ as the average number of DOFs in the VBA-results. \\
Table \ref{tab:ndofs_mu_uref} (a) clearly demonstrates the consistency of the scheme, whereby discrepancies as low as $\sim 0.5 \times 10^{-8}$ are achieved. Comparing the VBA to the SBA results, tables \ref{tab:ndofs_mu_uref} (a) and \ref{tab:minJ_nelems_uref} (a) demonstrate that VBA outperforms SBA in terms of accuracy, where an up to $\sim 10$-fold error reduction of VBA over SBA can be observed. The total number of iterations required until convergence is achieved is comparable for VBA and SBA and never exceeds the number of four iterations.

\begin{table}[h!]
\centering
\subfloat[$\#$DOFs for $\mathbf{x}^{\balpha}_h$.]{
\begin{tabular}{l|c|c|c}
  \backslashbox{\small $n_e$}{\small $u_{\text{ref}}$}  & $0$ & $1$ & $2$ \\ \hline
$12$ & $450$ & $450$ & $450$  \\ 
$16$ & $722$ & $722$ & $722$  \\
$23$ & $1352$ & $1352$ & $1352$
\end{tabular}} \qquad
\subfloat[$\#$DOFs for $u^{\balpha}_h$.]{
\begin{tabular}{l|c|c|c}
  \backslashbox{\small $n_e$}{\small $u_{\text{ref}}$}  & $0$ & $1$ & $2$ \\ \hline
$12$ & $225$ & $729$ & $2601$  \\ 
$16$ & $361$ & $1225$ & $4489$  \\
$23$ & $676$ & $2401$ & $9025$
\end{tabular} }
\caption{Tables showing the number of DOFs involved in computing $\mathbf{x}^{\balpha}_h$ (a) and $u^{\balpha}_h$ (b) in the SBA case for various combinations of $(n_e, u_{\text{ref}})$.}
\label{tab:ndofs_nelems_uref}
\end{table}

\subsection{Designing a Cooling Element}
\begin{figure}[h!]
\centering
\includegraphics[width=.8 \linewidth]{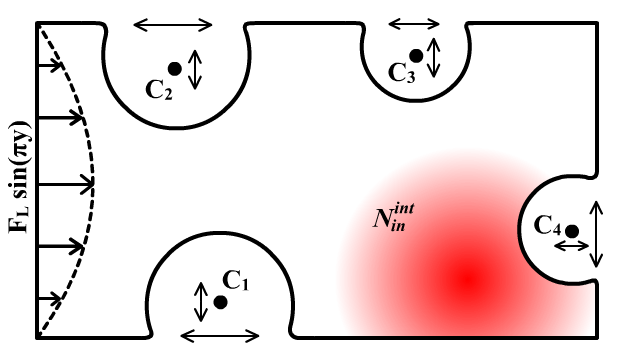}
\caption{The cooling element design template. Here, the centers of the active coolers are depicted by small black dots. Their positions constitute degrees of freedom in the design space, as well as the radii.}
\label{fig:cooling_element}
\end{figure}

\begin{figure}[h!]
  \centering
\includegraphics[width=.9\linewidth]{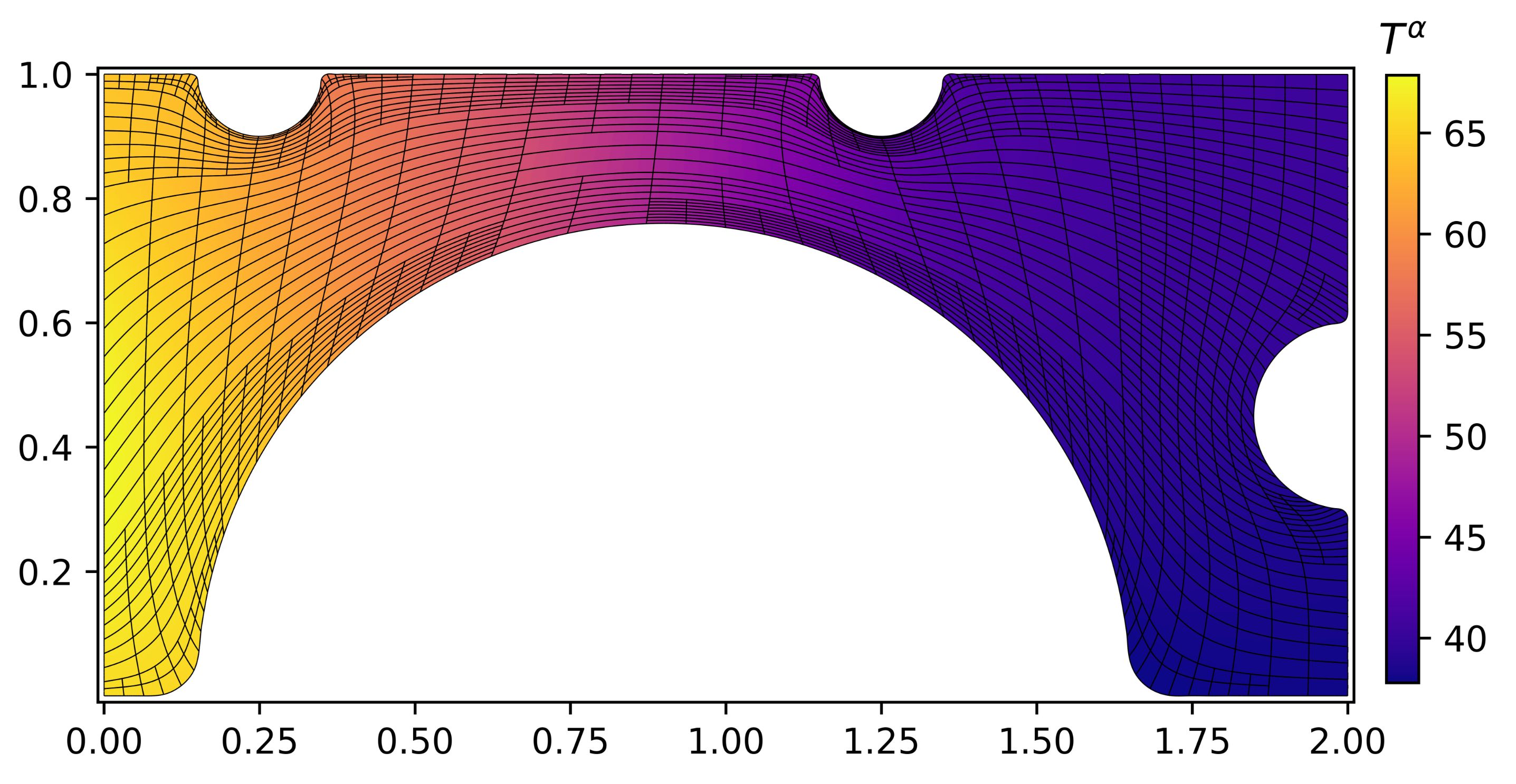}
\caption{The initial guess passed to the minimization routine.}
\label{fig:cooling_initial}
\end{figure}

We are considering the design of a cooling element of dimension $\Delta x = 2$ and $\Delta y = 1$. In this example, there are four active coolers whose positions can slide in the direction tangential to $\partial \Omega^{\balpha}$ and to a lesser extend in the normal direction (see Figure \ref{fig:cooling_element}). Further degrees of freedom are their radii $R_i$. Hence, the state vector is given by $\balpha = (\mathbf{x}_1, \mathbf{x}_2, \mathbf{x}_3, \mathbf{x}_4, R_1, R_2, R_3, R_4)$, which is comprised of $12$ DOFs. The surface cooling rate for the $i$-th active cooler $C_i$ reads:
\begin{align}
\label{eq:h_i}
    h^i_{-}(\mathbf{x}) = \frac{1}{20} \frac{R_i^3}{ \left\| \mathbf{x} - \mathbf{x}_i \right \|^2} \left( u^{\balpha}_h - T_{\infty} \right), \quad \text{where} \quad T_{\infty} = 0 \quad \text{denotes the ambient temperature}.
\end{align}
A heat source delivers a constant heat influx given by
\begin{align}
    N_{\text{tot}} = N_{\text{in}}^L + N_{\text{in}}^{\text{int}},
\end{align}
where $N_{\text{in}}^L$ denotes the influx at the left boundary, while $N_{\text{in}}^{\text{int}}$ denotes the influx delivered directly to the cooling element through an additional source term which satisfies
\begin{align}
\label{eq:Nin_int}
    N_{\text{in}}^{\text{int}} = \int \limits_{\Omega^{\balpha}} A \exp{\left(- \frac{\left \| \mathbf{x} - \mathbf{x}_0 \right\|^2}{2 \sigma^2} \right)} \mathrm{d} S, \quad \text{where} \quad A = \frac{N_{\text{tot}}}{4 \pi \sigma^2}, \quad \mathbf{x}_0 = (1.5, 0.25)^T \quad \text{and} \quad \mu = 0.1.
\end{align}
Note that changing the domain of integration from $\Omega^{\balpha}$ to $\mathbb{R}^2$ in (\ref{eq:Nin_int}) yields a value of $N_{\text{tot}} / 2$. As $N_{\text{tot}}$ is a constant quantity, we necessarily have $N_{\text{in}}^L = N_{\text{tot}} - N_{\text{in}}^{\text{int}}$. The surface heat flux density $h_L:\gamma^{\balpha}_L \rightarrow \mathbb{R}$ at the left boundary $\gamma^{\balpha}_L$ is of the form $h_L(x_2) = F_L(\Omega^{\balpha}) \sin{(\pi x_2)}$. Therefore, we have
\begin{align}
\label{eq:F_L}
    N_{\text{in}}^L = N_{\text{tot}} - N_{\text{in}}^{\text{int}} = \int \limits_{\gamma^{\balpha}_L} F_L \sin{(\pi x_2)} \mathrm{d} \gamma \implies F_L(\Omega^{\balpha}) = \frac{N_{\text{tot}} \pi}{2} \left( 1 - \int \limits_{\Omega^{\balpha}} \frac{1}{4 \pi \sigma^2} \exp{\left(- \frac{\left \| \mathbf{x} - \mathbf{x}_0 \right\|^2}{2 \sigma^2} \right)} \mathrm{d} S \right).
\end{align}
The relationship between $u^{\balpha}$ and the (uniform) temperature of the heat source $T^{\balpha}$ reads:
\begin{align}
\label{eq:Ntot_of_T_u}
    N_{\text{tot}} = A_1(\Omega^{\balpha}) \int \limits_{\gamma^{\balpha}_L} (T^{\balpha} - u^{\balpha}) \sin{(\pi x_2)} \mathrm{d} \gamma + A_2(\Omega^{\balpha}) \int \limits_{\Omega^{\balpha}} (T^{\balpha} - u^{\balpha}) \exp{\left(- \frac{\left \| \mathbf{x} - \mathbf{x}_0 \right\|^2}{2 \sigma^2} \right)} \mathrm{d} S,
\end{align}
where
\begin{align}
    A_1(\Omega^{\balpha}) = \frac{\pi}{2} \left( 1 - \frac{W(\Omega^{\balpha})}{4 \pi \sigma^2} \right) \quad \text{and} \quad A_2(\Omega^{\balpha}) = \frac{W(\Omega^{\balpha})}{8 \pi^2 \sigma^4}, \quad \text{with} \quad W(\Omega^{\balpha}) = \int \limits_{\Omega^{\balpha}} \exp{\left(- \frac{\left \| \mathbf{x} - \mathbf{x}_0 \right\|^2}{2 \sigma^2} \right)} \mathrm{d} S.
\end{align}
Inverting (\ref{eq:Ntot_of_T_u}) gives:
\begin{align}
\label{eq:T_of_u}
    T^{\balpha}(u^{\balpha}, \Omega^{\balpha}) = \frac{N_{\text{tot}} + A_1(\Omega^{\balpha}) \int_{\gamma^{\balpha}_L} u^{\balpha} \sin{(\pi x_2)} \mathrm{d} \gamma + A_2(\Omega^{\balpha}) \int_{\Omega^{\balpha}} u^{\balpha} \exp{\left(- \frac{\left \| \mathbf{x} - \mathbf{x}_0 \right\|^2}{2 \sigma^2} \right)} \mathrm{d} S}{\frac{2}{\pi} A_1(\Omega^{\balpha}) + W(\Omega^{\balpha}) A_2(\Omega^{\balpha})}.
\end{align}
\begin{remark}
The rationale behind $A_1(\Omega^{\balpha})$ and $A_2(\Omega^{\balpha})$ in (\ref{eq:T_of_u}) can be understood as follows: Given $N_{\text{tot}} = 1$, suppose the cooling element width were to be contracted to $\Delta x \rightarrow 0$. We have $\lim_{\Delta x \rightarrow 0} W(\Omega^{\balpha}) = 0$. In the limit, the temperature of the heat source should be fully determined by the first term on the right hand side of (\ref{eq:Ntot_of_T_u}), which is the case because $\lim_{\Delta x \rightarrow 0} (A_1, A_2) = (\tfrac{\pi}{2}, 0)$. As such, a constant influx of $1 = N_{\text{tot}} = N_{\text{in}}^L$ means $T^{\balpha} - u^{\balpha}\vert_{\gamma^{\balpha}_L} = 1$. Conversely, suppose $\Omega^{\balpha}$ were to be replaced by $\mathbb{R}^2$. Then, the dependency is divided equally among both terms since with $W(\Omega^{\balpha}) = 2 \pi \sigma^2$, $(A_1, A_2)= (\tfrac{\pi}{4}, \frac{1}{4 \pi \sigma^2})$. So, for $T^{\balpha} - u^{\balpha} = 1$, both terms contribute the same factor of $\tfrac{1}{2}$ to the right hand side of (\ref{eq:Ntot_of_T_u}).
\end{remark}
The weak state equation is based on the following PDE-problem:
\begin{align}
    - d \Delta u^{\balpha} &= -f u^{\balpha} + A \exp{\left(- \frac{\left \| \mathbf{x} - \mathbf{x}_0 \right\|^2}{2 \sigma^2} \right)} \nonumber \\
    \text{s.t.} \quad d \left. \frac{\partial u^{\balpha}}{\partial \mathbf{n}} \right \vert_{\partial \Omega^{\balpha}} & = \left\{
        \begin{array}{ll}
            - h_{\text{cooling}} + F_L \sin(\pi x_2) & \quad \mathbf{x} \in \bar{\gamma}_{L}^{\balpha} \\
            - h_{\text{cooling}} & \quad \mathbf{x} \in \partial \Omega^{\balpha} \setminus \bar{\gamma}_{L}^{\balpha}
        \end{array}
    \right.,
\end{align}
where
\begin{align}
    h_{\text{cooling}} = \sum_{i=1}^4 h^{i}_{-} \quad \text{and} \quad f = 10^{-3} \quad \text{denotes the internal dissipation rate}.
\end{align}
The $i$-th entry of the discretized weak state equation reads:
\begin{align}
    \left( \mathbf{B}^{\balpha}_h \right)_i = d \left( \nabla u^{\balpha}_h, \nabla \phi_i \right)_{\Omega^{\balpha}_h} - \int \limits_{\Omega^{\balpha}_h} f u^{\balpha} \phi_i \mathrm{d} S +  \int \limits_{\Omega^{\balpha}_h} A \exp{\left(- \frac{\left \| \mathbf{x} - \mathbf{x}_0 \right\|^2}{2 \sigma^2} \right) \left( \frac{\pi}{2} \tilde{\phi}_i  - \phi_i \right)} \mathrm{d}S + \sum_{j=1}^4 \int \limits_{\partial \Omega^{\balpha}_h} \phi_i h_{-}^j \mathrm{d} \gamma - \frac{\pi}{2} N_{\text{tot}} \tilde{\phi}_i,
\end{align}
with $A$ as in (\ref{eq:Nin_int}), $d = 0.8$ and
\begin{align}
    \tilde{\phi}_i = \int \limits_{\gamma_L^{\balpha}} \phi_i \sin(\pi y) \mathrm{d} \gamma.
\end{align}
We are minimizing the manufacturing costs of the cooling element such that the heat source temperature does not exceed the value of $T_{\text{max}} = 80$. The problem reads:
\begin{align}
\label{eq:Cooling_Element_Optimization}
\begin{split}
    J\left(u^{\balpha}, \Omega^{\balpha}, \balpha \right) & \rightarrow \min_{\balpha} \\
    \text{s.t.} \quad T_{\text{max}} - T^{\balpha} & \geq 0 \\
    \balpha & \in \boldsymbol{\mathcal{\lambda}},
\end{split}
\end{align}
where
\begin{align}
    J\left(u^{\balpha}, \Omega^{\balpha}, \balpha \right) = \int \limits_{\Omega^{\balpha}} 1 \mathrm{d}S + \sum \limits_{i=1}^4 C_{\text{CE}} R_i^2, \quad \text{with} \quad C_{\text{CE}} = \frac{100}{\pi}.
\end{align}
Furthermore, the feasible design space $\boldsymbol{\mathcal{\lambda}}$ is the space of all $\balpha$ such that the active coolers do not overlap and the genus of $\Omega^{\balpha}$ does not change (allowing for shape optimization without topology changes). This leads to a total of $30$ (partly nonlinear) inequalities. \\
A major challenge is deciding where to place the active coolers and what radii to use. Increasing the radius means additional cooling but also additional manufacturing costs and decreased cooling element area, decreasing the heat capacity and the channel heat conductivity. Furthermore, placing a cooler close to the internal heat source (see (\ref{eq:Nin_int})) reduces the amount of internal influx, increasing the influx amplitude $F_L$ (see (\ref{eq:F_L})) at $\gamma^{\balpha}_L$ for compensation. \\
We are considering the case $N_{\text{in}} = 10$ and follow the same approach as in Section \ref{sect:Validation} with $\mu = 0.5 \times 10^{-3}$ and $u_{\text{ref}} = 1$. Since $\mathbf{x}^{\balpha}$ is a continuous function of the input state vector, we improve the efficiency by storing the tuples $(\balpha^i, \cA^i, \boldsymbol{\mathcal{V}}^{\balpha, i}_h)$ (see Section \ref{sect:Computational_Approach}) after each iteration. Whenever some $\balpha^i$ with $\| \balpha - \balpha^i \| < \epsilon$ is found in the database, the corresponding mapping $\mathbf{x}^{\balpha, i}_h \in \boldsymbol{\mathcal{V}}^{\balpha, i}_h$ is prolonged to the coarsest element segmentation of $\hat{\Omega}$ that is compatible with both $\boldsymbol{\mathcal{V}}^{\balpha, i}_h$ and the current $\boldsymbol{\mathcal{V}}^{\balpha}_h$. Upon completion, it is restricted to $\boldsymbol{\mathcal{V}}^{\balpha}_h$, which yields the vector $\cA^R = \left( \cB^R, \cI^R \right)^T$. The weights corresponding to the inner DOFs, $\cI^R$, are extracted and then used as an initial guess for the root-finding problem (\ref{eq:x_weak_discretized}). We have noticed this to lead to a tremendous speedup, in particular during the last iterations, in which $\balpha$ varies only slightly. Hereby, the required number of iterations is reduced from typically four to as few as one.
\begin{remark}
This principle may be extended to higher than zeroth-order database interpolation.
\end{remark}
Here, we use $\epsilon = 0.05$. A feasible initial guess is created by picking one of the coolers and increasing its radius until $T^{\balpha} < T_{\text{max}}$. The initial design is depicted in Figure \ref{fig:cooling_initial}. As in section \ref{sect:Validation}, the routine that computes $J^{\balpha}_h$, $\mathrm{d}_{\balpha} J^{\balpha}_h$ and the constraints is passed to an SLSQP optimizer.

\begin{figure}[h!]
  \centering
  \subfloat[]{\includegraphics[width=.45\linewidth]{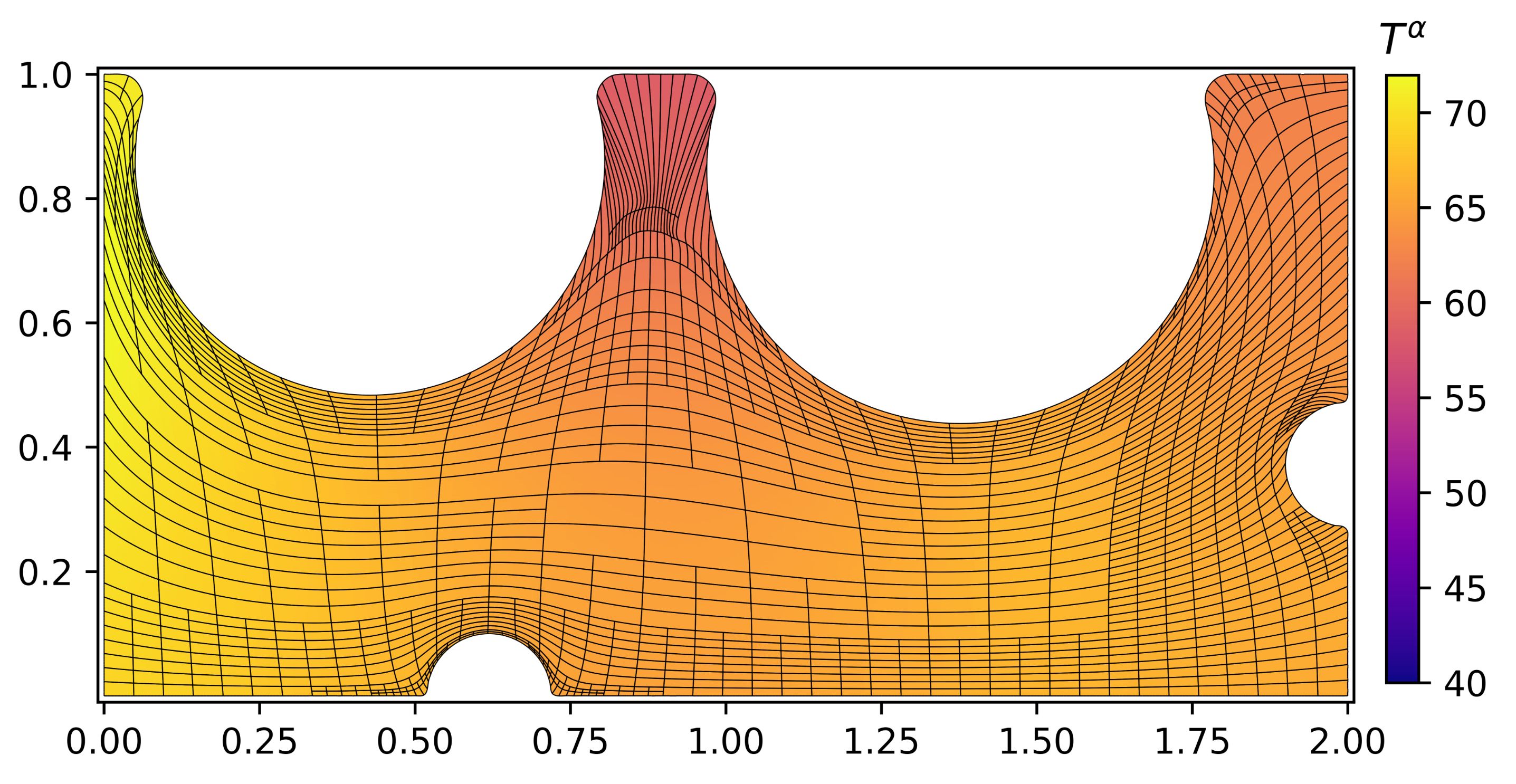}\label{fig:cooling_6}} $\quad$
  \subfloat[]{\includegraphics[width=.45\linewidth]{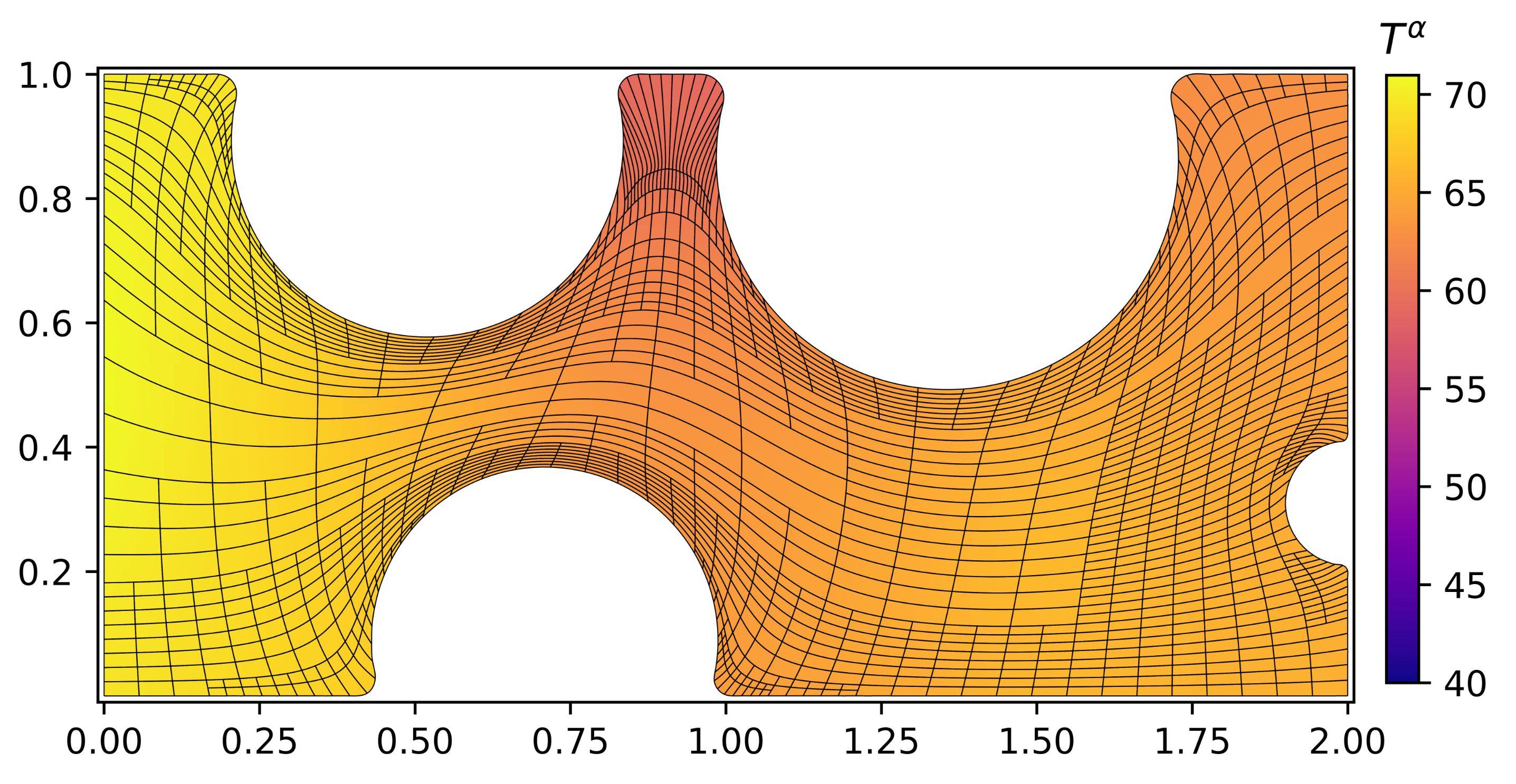}\label{fig:cooling_12}} \\
  \subfloat[]{\includegraphics[width=.45\linewidth]{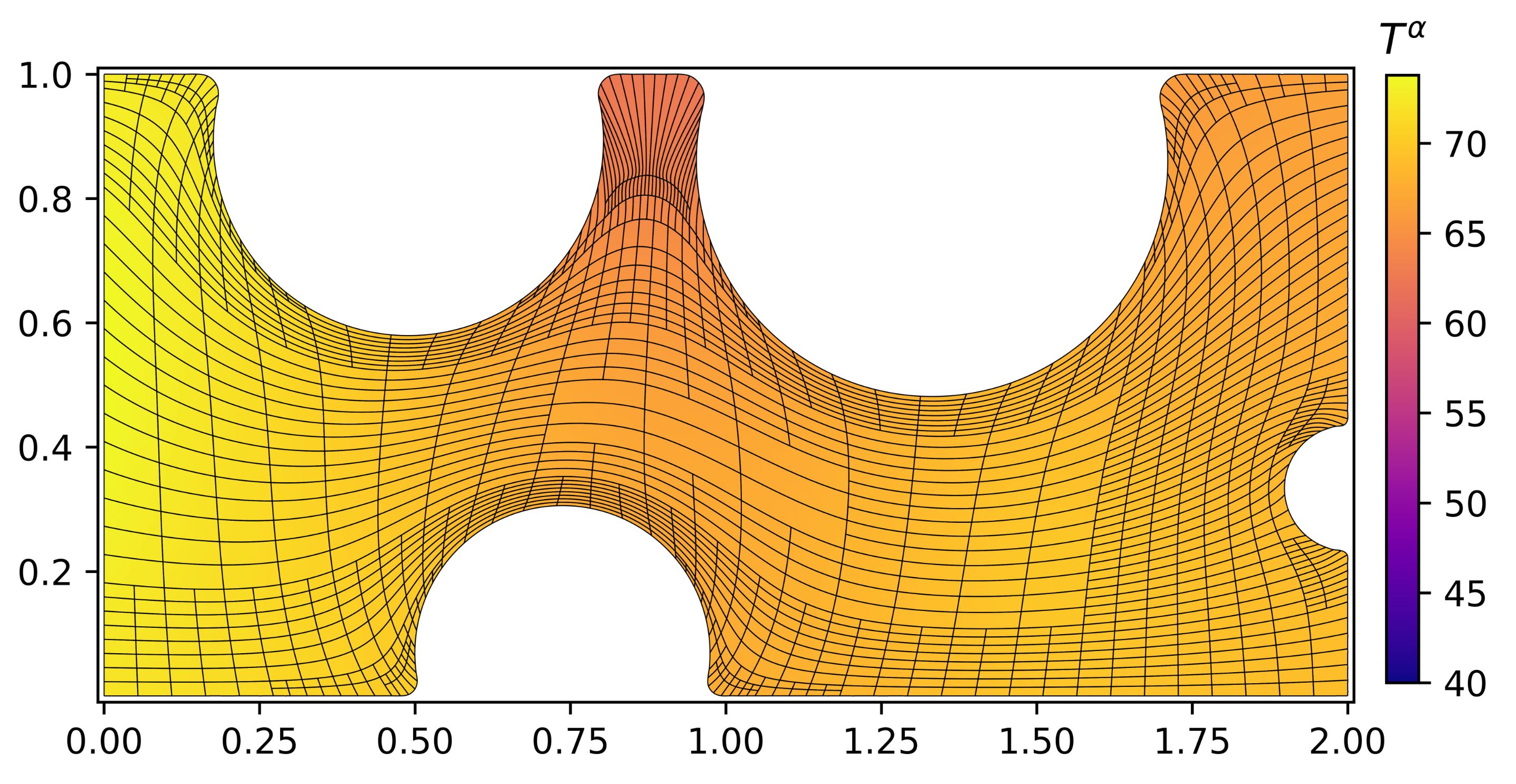}\label{fig:cooling_18}} $\quad$
  \subfloat[]{\includegraphics[width=.45\linewidth]{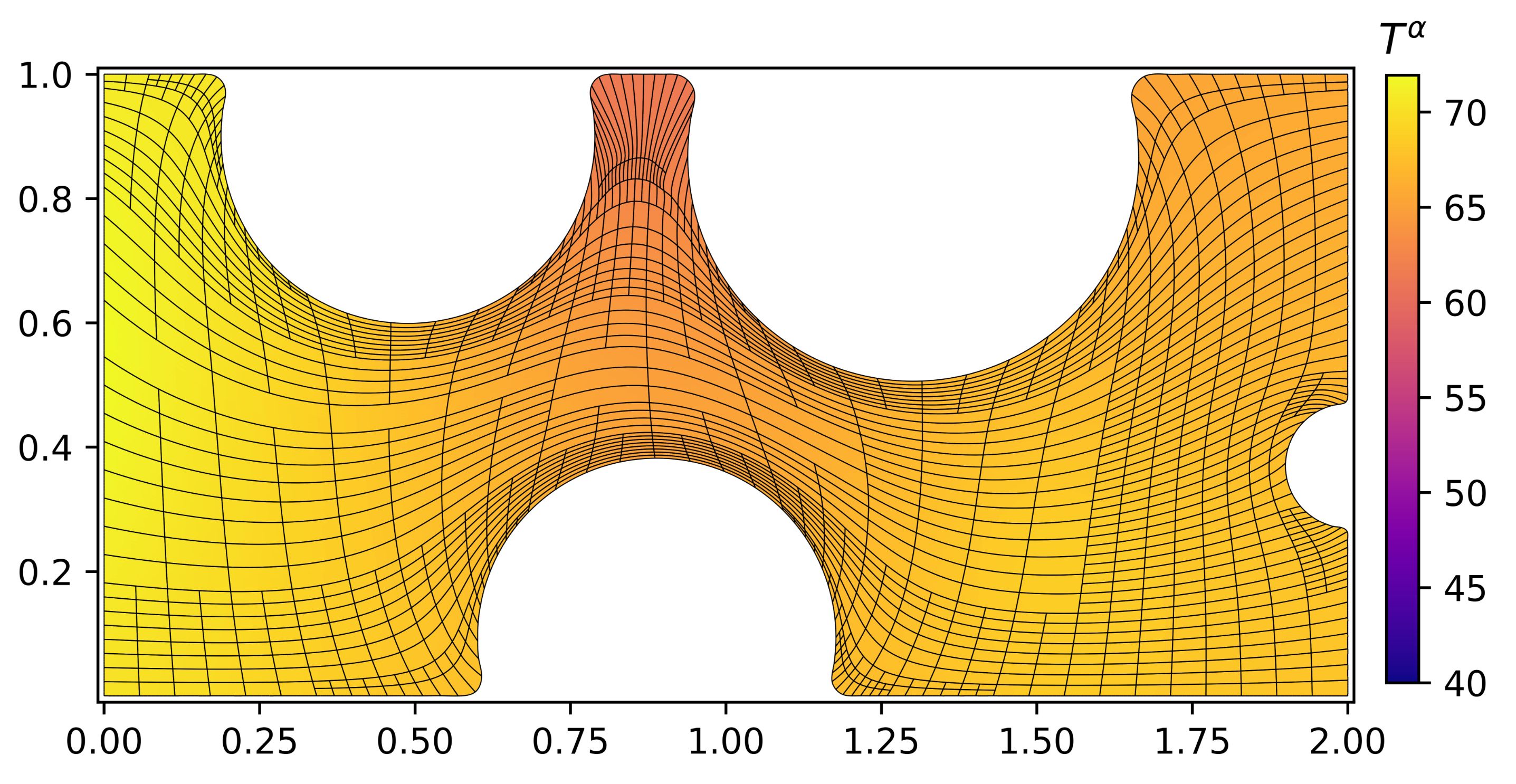}\label{fig:cooling_24}}
  \caption{The cooling element after $4$ (a), $7$ (b), $10$ (c) and $13$ (d) iterations.}
  \label{fig:cooling_6_24}
\end{figure}

\begin{figure}[h!]
  \centering
\includegraphics[width=.9\linewidth]{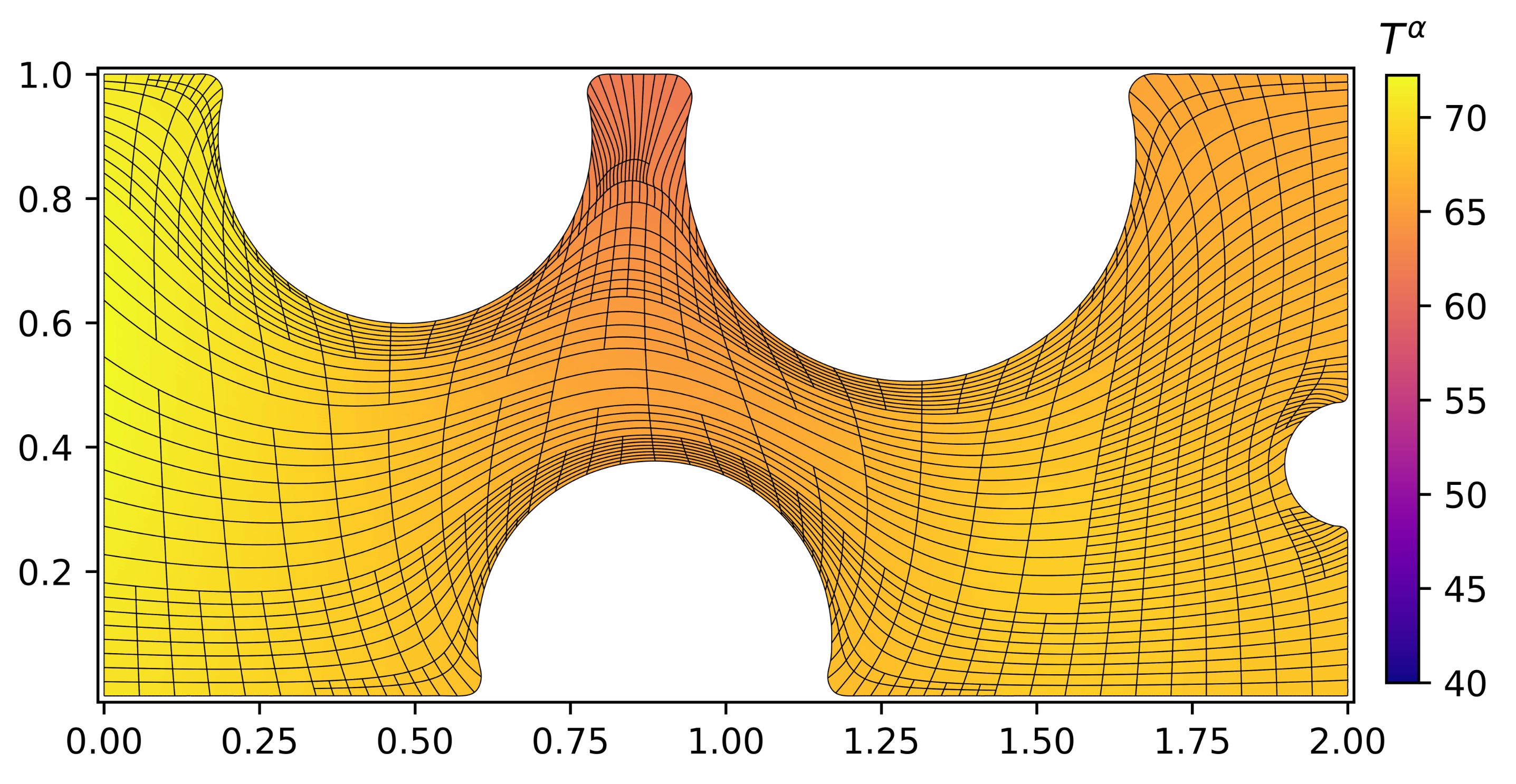}
\caption{The final cooling element design after $15$ iterations.}
\label{fig:cooling_final}
\end{figure}

Figures \ref{fig:cooling_6_24} (a) to \ref{fig:cooling_6_24} (d) show the cooling element after $4$, $7$, $10$ and $13$ iterations. Convergence is reached after $15$ iterations and the corresponding design is depicted in Figure \ref{fig:cooling_final}. The final design reduces the manufacturing costs from the initial $J^{\balpha}_h = 10.66$ to $J^{\balpha}_h = 6.29$. \\
A striking difference between the initial and all intermediate designs is the improved heat conductivity within the channel, leading to a more homogeneous temperature (and one that is higher on average). This is not surprising. As the cooling efficiency is linear in the difference between the temperature at the boundaries and the ambient temperature $T_{\infty} = 0$, a higher average temperature implies higher average cooling efficiency. \\
The final design places a modestly-sized cooler $\mathcal{C}_1$ at the center of the southern boundary and a similarly-sized cooler $\mathcal{C}_2$ at the western part of the northern boundary. To its right, a slightly larger cooler $\mathcal{C}_3$ is placed while a small cooler $\mathcal{C}_4$ is placed at eastern boundary close to $x_2 = 0.4$.

Compared to the initial design, one big cooler has been replaced by several modestly-sized ones, improving the channel heat conductivity and by that the cooling cost efficiency. The slightly larger size of $\mathcal{C}_3$ compared to $\mathcal{C}_2$ can be explained by the internal heat source centered at $\mathbf{x}_0 = (1.5, 0.25)^T$. The small radius of $\mathcal{C}_4$ may be explained by the fact that increasing its size reduces the amount of internal influx area, leading to a larger influx at $\gamma_L^{\balpha}$ instead. As such, we regard the final design as plausible, adding more credibility to the proposed numerical scheme.
\section{Conclusion}
In this manuscript, we proposed an IGA-based shape optimization algorithm in which the parameterization is included in the problem formulation in the form of an additional PDE-constraint. This has enabled us to derive a fully symbolical expression for the gradient of the objective function, allowing for gradient-based optimization. The discretization of the equations has been accomplished with the so-called \textit{variable basis approach} (VBA) in which a new THB-spline basis is chosen during each iteration based on the current requirements, such as accurately resolving the geometry contours and particular features of the state equation solution. This leads to a highly flexible scheme in which folding due to numerical truncation is automatically repaired through THB-enabled local refinement. \\
We have tested the scheme by applying it to two examples. In the first example, we compared the numerical solution to the known exact solution and concluded that the scheme is consistent. Comparing the VBA-approach to an approach in which the basis is taken static (SBA) furthermore revealed that VBA-enabled feature-based refinement leads to a $\sim 10$-fold error reduction over SBA at a comparable total number of DOFs. This discrepancy may be further increased by employing more proficient a priori and a posteriori refinement techniques. In the second example, we considered the design of a cooling element. Unlike in the first example, the exact minimizer was unknown, however, the optimization routine converged to a design that we consider plausible. In both cases, the scheme succeeded in fully automatically parameterizing a wide range of geometries which would be too complex for other symbolically-differentiable parameterization strategies (such as Coon's Patch) at the expense of leading to a nonlinear problem. \\
Finally, we briefly discussed possible memory-saving strategies for large-scale optimization and (possible) future implementations of the scheme with support for volumetric applications. Furthermore, the scheme is straightforwardly enhanced to support multipatch parameterizations by adopting the mixed FEM EGG algorithm introduced in \cite{hinz2019iga}.

\section*{Acknowledgements}
The authors gratefully acknowledge the research funding which was partly provided by the MOTOR project that has received funding from the European Unions Horizon 2020 research and innovation program under grant agreement No 678727. 

The research visit of Mr. A. Jaeschke at Delft University of Technology that allowed this cooperation was funded by The Polish National Agency for Academic Exchange under grant agreement No PPN/IWA/2018/1/00032.

\bibliographystyle{elsarticle-num}
\bibliography{bibliography}


\end{document}